\documentclass[10pt]{article}

\usepackage{amsfonts}
\usepackage{amsmath}
\usepackage{amssymb}
\usepackage{amsthm}

\usepackage{a4wide}
\usepackage{longtable}
\usepackage{graphicx}
\usepackage{multirow}

\usepackage{stmaryrd}

\usepackage[latin1]{inputenc}
\usepackage[english]{babel}
\usepackage[T1]{fontenc}
\usepackage{amstext,latexsym}

\usepackage{pgf,tikz}
\usepackage{mathrsfs}
\usetikzlibrary{arrows}

\newtheorem{theorem}{Theorem}[section]

\theoremstyle{definition}

\theoremstyle{remark}
\newtheorem{remark}{Remark}[section]

\usepackage[all]{xy}
\input xy
\xyoption {all}
\usepackage[pdftex,             %
    bookmarks=true,             % Signets
    hypertexnames=false,        % Corrige les warnings (?)
    bookmarksnumbered=true,     % Signets numérotés
    pdfpagemode=None,           % Signets/vignettes fermé à l'ouverture
    pdfstartview=FitH,          % La page prend toute la largeur
    pdfpagelayout=SinglePage,   % Vue par page
    colorlinks=true,            % Liens en couleur
    urlcolor=blue,              % Couleur des liens externes
    citecolor=blue,             % Couleur des num biblio
    pdfborder={0 0 0}           % Style de bordure
    ]{hyperref}

\begin{document}

\title{Structure of relative genus fields \\ of cubic Kummer extensions}

\author{S. AOUISSI\footnote{Corresponding author. E.mail: aouissi.siham@gmail.com}, A. AZIZI, M. C. ISMAILI, D. C. MAYER and M. TALBI}
\date{Friday, 13 October 2023}
\maketitle

%-----------------------------------------------------------------------

\medskip
\noindent
{\bf Abstract:}
Let $N=K(\sqrt[3]{D})$ be a cubic Kummer extension
of the cyclotomic field $K=\mathbb{Q}(\zeta_3)$,
containing a primitive cube root of unity $\zeta_3$,
with cube free integer radicand $D>1$.
Denote by $f$ the conductor of the abelian extension $N/K$, 
and by $N^{\ast}$ the relative genus field of $N/K$.
The aim of the present work is
to find out all positive integers $D$ and conductors $f$ such that
the genus group
$\operatorname{Gal}\left(N^{\ast}/N\right)\cong
\mathbb{Z}/3\mathbb{Z}\times\mathbb{Z}/3\mathbb{Z}$
is elementary bicyclic.

\bigskip
\noindent
{\bf Keywords:} {Pure cubic field, cubic Kummer extension,
relative genus field, group of ambiguous ideal classes,
$3$-rank, primitive ambiguous principal ideals, multiplicity of conductors.}

\bigskip
\noindent
{\bf Mathematics Subject Classification 2010:} {11R11, 11R16, 11R20, 11R27, 11R29, 11R37.}
\bigskip

%%%%%%%%%%%%%%%%%%%%%%%%%%%%%%%%%%%%%%%%%%%%%%%%%%%%%%%%%%%%%%%%%%%%%%%%%%%%%%%%%%%%%%

\section{Introduction}
\label{s:Intro}
Let
$K=\mathbb{Q}(\zeta_3)$, where $\zeta_3$ is a primitive cube root of unity,
$N=\mathbb{Q}(\sqrt[3]{D},\zeta_3)$, where $D>1$ is a cube free integer radicand,
$\operatorname{Gal}(N/K)=\langle\sigma\rangle$ the cyclic relative Galois group of $N/K$,
$Cl_{3}(N)$ be the \(3\)-class group of \(N\),
$N_3^{(1)}$ be the maximal abelian unramified $3$-extension of $N$,
$f$ the conductor of the abelian extension $N/K$ and $m$ its multiplicity.
Let $N^{\ast}$ be the maximal abelian extension of $K$ contained in $N_3^{(1)}$,
which is called the \textit{relative genus field} of $N/K$
(see \cite{Ge1975}, \cite{Ge1976} or \cite{Hz}).
\paragraph{}
One tries to determine the unramified $3$-sub-extensions of $N_3^{(1)}/N$
and then, according to class field theory,
extract information about $Cl_{3}(N)$, its rank,
and the $3$-class field tower of $N$.
One way to do this is by asking for the structure of the relative genus field $N^{\ast}$ of $N/K$.
In our recent work \cite{AMITA}, we implemented Gerth's methods \cite{Ge1975} and \cite{Ge1976}
to determine the rank of the group $Cl_{3}^{(\sigma)}(N)$ of ambiguous ideal classes of $N/K$ and
obtained all integers $D>1$ for which
$\operatorname{Gal}(N^{\ast}/N)\cong\mathbb{Z}/3\mathbb{Z}$ is cyclic of order $3$.
The purpose of the present work is to find all positive integers $D$ and conductors $f$ such that
$\operatorname{Gal}(N^{\ast}/N)\cong\mathbb{Z}/3\mathbb{Z}\times\mathbb{Z}/3\mathbb{Z}$
is elementary bicyclic.
In fact, we shall prove the following Main Theorem:

%-----------------------------------------------------------------------

\begin{theorem}
\label{thm:Rank2}
Let $N=\mathbb{Q}(\sqrt[3]{D},\zeta_3)$, where $D$ is a cube free positive integer,
$K= \mathbb{Q}(\zeta_3)$, and  $N^{\ast}$ the relative genus field of $N/K$.
Then $\operatorname{Gal}(N^{\ast}/N)\cong\mathbb{Z}/3\mathbb{Z}\times\mathbb{Z}/3\mathbb{Z}$
if and only if the integer $D$ can be written in one of the following thirteen forms:
\newpage
\begin{equation}
\label{eqn:Rank2}
D = \left\{
  \begin{array}{l l}
    %%%%%%%%% 1 %%%%%%%%%
    3^{e}p^{e_1} & \text{ with } p\equiv 1\,(\mathrm{mod}\,9), \\
    %%%%%%%%% 2 %%%%%%%%%
    p^{e_1}q^{f_1} & \text{ with } p\equiv -q\equiv 1\,(\mathrm{mod}\,9), \\
    %%%%%%%%% 3 %%%%%%%%%
    p^{e_1}q^{f_1} \not\equiv\pm 1\,(\mathrm{mod}\,9) & \text{ with } p\not\equiv 1\,(\mathrm{mod}\,9) \text{ or } q\not\equiv -1\,(\mathrm{mod}\,9), \\
    %%%%%%%%% 4 %%%%%%%%%
    3^{e}p^{e_1}q^{f_1} & \text{ with } p\not\equiv 1\,(\mathrm{mod}\,9) \text{ or } q\not\equiv -1\,(\mathrm{mod}\,9), \\
    %%%%%%%%% 5 %%%%%%%%%
    p^{e_1}q_1^{f_1}q_2^{f_2} \equiv\pm 1\,(\mathrm{mod}\,9) & \text{ with } p\equiv 1\,(\mathrm{mod}\,9) \text{ and } q_1,q_2\equiv 2 \text{ or } 5\,(\mathrm{mod}\,9), \\
    %%%%%%%%% 6 %%%%%%%%%
    p^{e_1}q_1^{f_1}q_2^{f_2} \equiv\pm 1\,(\mathrm{mod}\,9) & \text{ with } p,-q_1,-q_2\equiv 4 \text{ or } 7\,(\mathrm{mod}\,9), \\
    %%%%%%%%% 7 %%%%%%%%%
    p^{e_1}q_1^{f_1}q_2^{f_2} \equiv\pm 1\,(\mathrm{mod}\,9) & \text{ with } p,-q_2 \equiv 4 \text{ or } 7\,(\mathrm{mod}\,9) \text{ and } q_1 \equiv -1\,(\mathrm{mod}\,9), \\
    %%%%%%%%% 8 %%%%%%%%%
    p_1^{e_1}p_2^{e_2} \equiv\pm 1\,(\mathrm{mod}\,9) & \text{ such that } \exists i\in\lbrace 1,2\rbrace\mid p_{i}\not\equiv 1\,(\mathrm{mod}\,9), \\
    %%%%%%%%% 9 %%%%%%%%%
    3^{e}q_1^{f_1}q_2^{f_2} & \text{ with } q_1\equiv q_2\equiv -1\,(\mathrm{mod}\,9), \\
    %%%%%%%% 10 %%%%%%%%%
    q_1^{f_1}q_2^{f_2}q_3^{f_3} & \text{ with } q_1\equiv q_2\equiv q_3\equiv -1\,(\mathrm{mod}\,9), \\
    %%%%%%%% 11 %%%%%%%%%
    q_1^{f_1}q_2^{f_2}q_3^{f_3} \not\equiv\pm 1\,(\mathrm{mod}\,9) & \text{ such that } \exists i\in\lbrace 1,2,3\rbrace\mid q_{i}\not\equiv -1\,(\mathrm{mod}\,9), \\
    %%%%%%%% 12 %%%%%%%%%
    3^{e}q_1^{f_1}q_2^{f_2}q_3^{f_3} & \text{ such that } \exists i\in\lbrace 1,2,3\rbrace\mid q_{i}\not\equiv -1\,(\mathrm{mod}\,9), \\
    %%%%%%%% 13 %%%%%%%%%
    q_1^{f_1}q_2^{f_2}q_3^{f_3}q_4^{f_4} \equiv\pm 1\,(\mathrm{mod}\,9) & \text{ such that } \exists i\in\lbrace 1,2,3,4\rbrace\mid q_{i}\not\equiv -1\,(\mathrm{mod}\,9),
  \end{array} \right.
\end{equation}
where $e, e_1, e_2, f_1, f_2, f_3$ and $f_4$ are positive integers equal to $1$ or $2$.
\end{theorem}
 
%-----------------------------------------------------------------------

\paragraph{}
As opposed to the radicands $D$ of the shape in \cite[Thm. 1.1, p 251]{AMITA},
the possible prime factorizations of $D$ in our Main Theorem \ref{thm:Rank2}
are more complicated. For background see section \ref{s:Background}.
 
\paragraph{}
In section \ref{s:Context},
we give Theorem \ref{thm:Ismaili1} and Theorem \ref{thm:Ismaili2} which show that
pure cubic fields $L=\mathbb{Q}(\sqrt[3]{D})$ can be collected in multiplets \((L_1,\ldots,L_m)\)
sharing a common conductor $f$ with multiplicity $m$
and a common type of ambiguous class group $Cl_{3}^{(\sigma)}(N)$ of $N$.
At the beginning of section \ref{s:Proof}, where Theorem \ref{thm:Rank2} is proved,
we restrict ourselves to those results that will be needed in this paper. 
More information on $3$-class groups can be found in \cite{Ge1973}, \cite{Ge1975}, \cite{Ge1976}, \cite{Ky1973}, \cite{Ky1974}, \cite{Ky}, and \cite{Ky1977}.
For the prime decomposition in a pure cubic field $\mathbb{Q}(\sqrt[3]{D})$,
we refer the reader to the papers \cite{BC1970}, \cite{BC1971}, \cite{Dk}, and \cite{Mk},
and for the prime ideal factorization rules of $\mathbb{Q}(\zeta_3)$, we refer to \cite{IlRo}.
In section \ref{ss:Split},
we present numerical examples concerning the structure of the $3$-class groups of pure cubic fields $L$
and of their Galois closures $N$ in the case where conductors $f$ contain splitting prime divisors.

\begin{center}
\textbf{Notations:}
\end{center}
\begin{itemize}
 \item $p$ and $q,\ell$ are prime numbers such that $p\equiv 1 \pmod 3$ and  $q,\ell\equiv -1\pmod 3$;
 \item $L=\mathbb{Q}(\sqrt[3]{D})$ is a pure cubic field, where $D>1$ is a cube free positive integer;
 \item $K=\mathbb{Q}(\zeta_3)$, where $\zeta_3=e^{2i\pi/3}$ denotes a primitive third root of unity;
 \item $N=\mathbb{Q}(\sqrt[3]{D},\zeta_3)$ is the normal closure of $L$;
 \item $f$ is the conductor of the abelian extension $N/K$;
 \item $\langle\tau\rangle=\operatorname{Gal}(N/L)$ such that $\tau^2=id$, $\tau(\zeta_3)=\zeta_3^2$ and $\tau(\sqrt[3]{D})=\sqrt[3]{D}$;
 \item $\langle\sigma\rangle=\operatorname{Gal}(N/K)$ such that $\sigma^3=id$, $\sigma(\zeta_3)=\zeta_3$,
   $\sigma(\sqrt[3]{D})=\zeta_3\sqrt[3]{D}$ and $\tau\sigma=\sigma^2\tau$;
 \item $\lambda=1-\zeta_3$ and $\pi,\rho$ are prime elements of $K$;
 \item $q^{\ast}=1$ or \(0\) according to whether \(\zeta_{3}\) is or is not norm of an element of \(N\setminus\lbrace 0\rbrace\);
 \item \(t\) denotes the number of prime ideals ramified in \(N/K\).
 \item For an algebraic number field $F$:
  \begin{itemize}
  \item $\mathcal{O}_{F}$, $E_{F}$ : the ring of integers and the group of units of $F$;
  \item $Cl_{3}(F)$, $F_3^{(1)}$ : the $3$-class group and the Hilbert $3$-class field of $F$.
  \end{itemize}
\end{itemize}

%%%%%%%%%%%%%%%%%%%%%%%%%%%%%%%%%%%

\newpage
%-----------------------------------------------------------------------

\section{Background}
\label{s:Background}

\noindent
For the convenience of the reader,
we collect all formulas for pure cubic fields
$L=\mathbb{Q}(\sqrt[3]{D})$
and their normal closures
$N=\mathbb{Q}(\sqrt[3]{D},\zeta_3)$
on which subsequent proofs and exposition build up.
Generally, we adopt the notation of Gerth
\cite{Ge1975}.
In particular, the prime decomposition of the cube free integer \textit{radicand} $D$
is written in the following form (Equation $(3.2)$ of \cite[p. 55]{Ge1975}):
\begin{eqnarray}
\label{eqn:PrimeDecRad}
D &=& 3^e\cdot p_1^{e_1}\cdots p_v^{e_v}\cdot p_{v+1}^{e_{v+1}}\cdots p_w^{e_w}\cdot q_1^{f_1}\cdots q_I^{f_I}\cdot q_{I+1}^{f_{I+1}}\cdots q_J^{f_J},
\end{eqnarray}
where $p_i$ and $q_i$ are rational prime numbers such that:
\[
\left\{
\begin{array}{ll}
 p_i \equiv 1\,(\mathrm{mod}\,9), & \quad \text{ for } \quad 1\leq i\leq v, \\
 p_i \equiv 4 \text{ or } 7\,(\mathrm{mod}\,9), & \quad \text{ for } \quad v+1\leq i\leq w,\\
 q_i \equiv 8\,(\mathrm{mod}\,9), & \quad \text{ for } \quad 1\leq i\leq I,\\
 q_i \equiv 2 \text{ or } 5\,(\mathrm{mod}\,9), & \quad \text{ for } \quad I+1\leq i\leq J,\\
 e_i = 1 \text{ or } 2, & \quad \text{ for } \quad 1\leq i\leq w,\\
 f_i = 1 \text{ or } 2, & \quad \text{ for } \quad 1\leq i\leq J,\\
 e   = 0 \text{ or } 1 \text{ or } 2.
\end{array}
\right.
\]
According to
\cite[Eqn. (2.5)--(2.6), p. 255]{AMITA}
the corresponding prime decomposition of the class field theoretic \textit{conductor} \(f\) of the abelian relative extension \(N/K\)
\cite{Ha1930}
is given by
\begin{eqnarray}
\label{eqn:PrimeDecCond}
f &=& 3^\varepsilon\cdot p_1\cdots p_v\cdot p_{v+1}\cdots p_w\cdot q_1\cdots q_I\cdot q_{I+1}\cdots q_J.
\end{eqnarray}
Here the exponent \(\varepsilon=v_3(f)\) of the distinguished prime \(3\) characterizes the \textit{Dedekind species} of \(L\): 
\begin{equation}
\label{eqn:Species}
\varepsilon=
\begin{cases}
2 & \text{ if } e\in\lbrace 1,2\rbrace \text{ (Species 1a)}, \\
1 & \text{ if } e=0 \text{ and } D\not\equiv\pm 1\,(\mathrm{mod}\,9) \text{ (Species 1b)}, \\
0 & \text{ if } e=0 \text{ and } D\equiv\pm 1\,(\mathrm{mod}\,9) \text{ (Species 2)}. \\
\end{cases}
\end{equation}
Additionally, we use the following conventions
for conjugate prime elements in the cyclotomic field \(K=\mathbb{Q}(\zeta_3)\)
which divide \textit{splitting} rational prime numbers:
\(p=\rho\rho^\tau\), if \(p\equiv 1\,(\mathrm{mod}\,9)\), but
\(p=\pi\pi^\tau\), if \(p\equiv 4,7\,(\mathrm{mod}\,9)\).
For \textit{inert} rational prime numbers,
we replace \(q\) by \(\ell\), if \(q\equiv 8\,(\mathrm{mod}\,9)\),
but we keep \(q\), if \(q\equiv 2,5\,(\mathrm{mod}\,9)\).

In order to determine the \(3\)-rank \(r\) of the ambiguous \(3\)-class group $Cl_{3}^{(\sigma)}(N)$ of \(N/K\),
Gerth used the general formula of Hasse
\cite[Thm. 13]{Ha1927}
in \cite[Eqn. (2.1), pp. 86--87]{Ge1976}
and, adapted to the particular situation \(K=\mathbb{Q}(\zeta_3)\) in
\cite[Prop. 5.1, pp. 92--93]{Ge1976},
in the \textit{implicit} form \(r=t+q^\ast-2\), \(t=J+2w(+1)\).
However, we only need the following entirely \textit{explicit} version
\cite[Lem. 3.1, p. 55]{Ge1975}
for the normal closure \(N\) of a pure cubic field \(L\),
in dependence on the Dedekind species \(v_3(f)\):
\begin{equation}
\label{eqn:AmbiguousRank}
r=
\begin{cases}
J+2w   & \text{ for } (v=w \text{ and } I=J) \text{ and } v_3(f)>0 \text{ (Species 1)}, \\
J+2w-1 & \text{ for }  (v<w \text{ or } I<J) \text{ and } v_3(f)>0 \text{ (Species 1)}, \\
J+2w-1 & \text{ for } (v=w \text{ and } I=J) \text{ and } v_3(f)=0 \text{ (Species 2)}, \\
J+2w-2 & \text{ for }  (v<w \text{ or } I<J) \text{ and } v_3(f)=0 \text{ (Species 2)}.
\end{cases}
\end{equation} 
According to
\cite[Thm. 2.1, p. 833]{Ma1992} and \cite[Cor. 3.2, p. 2219]{Ma2014},
the \textit{multiplicity} \(m(f)\) of the conductor \(f\),
that is, the number of pairwise non-isomorphic pure cubic fields \(L=\mathbb{Q}(\sqrt[3]{D})\)
sharing the common discriminant \(d_L=-3\cdot f^2\),
is given in terms of all,
respectively free (good) and restrictive (bad), prime divisors of \(f\),
i.e. \(T=w+J\), respectively \(G=v+I\), \(B=T-G\),
in dependence on the Dedekind species \(v_3(f)\) by
\begin{equation}
\label{eqn:Multiplicity}
m=m(f)=
\begin{cases}
2^T,              & \text{ if } v_3(f)=2 \text{ (Species 1a)}, \\
2^G\cdot X_B,     & \text{ if } v_3(f)=1 \text{ (Species 1b)}, \\
2^G\cdot X_{B-1}, & \text{ if } v_3(f)=0 \text{ (Species 2)},
\end{cases}
\end{equation}
where the sequence \((X_k)_{k\ge -1}\), \(X_k=\frac{1}{3}\lbrack 2^k-(-1)^k\rbrack\) starts with \(\frac{1}{2},0,1,1,3,5,11\).
In particular, the multiplicity is \textit{zero},
either if \(B=0\) for a field of species 1b
or if \(B=1\) for a field of species 2.

Now we give an \textit{elementary simultaneous proof} of
Honda's Theorem
\cite[Thm., \S\ 1, p. 8]{Ho1971},
our Main Theorem in
\cite[Thm. 1.1, p. 251]{AMITA},
and our Main Theorem
\ref{thm:Rank2}
in the present article.
The proof consists of a \textit{systematic combinatorial construction}
of all possible conductors \(f\),
ordered firstly by increasing number \(t\) of ramified prime ideals of \(N/K\),
and, for fixed \(t\), secondly by decreasing number \(w\) of primes splitting in \(K\).
For each conductor \(f\), the rank \(r\) and the multiplicity \(m\)
are determined with the aid of our collection of formulas.
On principle, this could be done by a computer search up to a given maximal \(t\),
if all deterministic rules concerning
the \(3\)-rank \(r\) of ambiguous class groups
\eqref{eqn:AmbiguousRank}
and the multiplicity \(m\) of conductors \(f\)
\eqref{eqn:Multiplicity}
are implemented in the program script.
The \textit{principal factorization type} is abbreviated by PFT
\cite[Thm. 2.1, p. 254]{AMITA}.

%\newpage
%--------------------------------------------------------------------------------

\renewcommand{\arraystretch}{1.2}

\begin{table}[ht]
\caption{Arithmetic invariants for conductors of Dedekind species 1a}
\label{tbl:Dedekind1a}
\begin{center}
{\normalsize
\begin{tabular}{|c|c|c|c||l||c|c|c||c|c|c|c||c|c|}
\hline
 \(J\) & \(w\) & \(I\) & \(v\) & \(f\)             & \(t\) & \(q^\ast\) & \(r\) & \(G\) & \(B\) & \(T\) & \(m\) & PFT              & Reference \\
\hline
 \(0\) & \(0\) & \(0\) & \(0\) & \(3^2\)                & \(1\) & \(1\) & \(0\) & \(0\) & \(0\) & \(0\) & \(1\) & \(\gamma\)       & \cite[Thm. 2.3 (1)]{AMITA} \\
 \(1\) & \(0\) & \(0\) & \(0\) & \(3^2q\)               & \(2\) & \(0\) & \(0\) & \(0\) & \(1\) & \(1\) & \(2\) & \(\beta\)        & \cite[Thm. 2.3 (4)]{AMITA} \\
 \(1\) & \(0\) & \(1\) & \(0\) & \(3^2\ell\)            & \(2\) & \(1\) & \(1\) & \(1\) & \(0\) & \(1\) & \(2\) & \(\beta,\gamma\) & \cite[Thm. 2.5 (1)]{AMITA} \\
 \(0\) & \(1\) & \(0\) & \(0\) & \(3^2\pi\pi^\tau\)     & \(3\) & \(0\) & \(1\) & \(0\) & \(1\) & \(1\) & \(2\) & \(\alpha,\beta\) & \cite[Thm. 2.4 (3)]{AMITA} \\
 \(0\) & \(1\) & \(0\) & \(1\) & \(3^2\rho\rho^\tau\)   & \(3\) & \(1\) & \(2\) & \(1\) & \(0\) & \(1\) & \(2\) & \(\alpha,\beta,\gamma\) & Thm. 3.1 (1) \\
 \(2\) & \(0\) & \(0\) & \(0\) & \(3^2q_1q_2\)          & \(3\) & \(0\) & \(1\) & \(0\) & \(2\) & \(2\) & \(4\) & \(\beta\)        & \cite[Thm. 2.5 (5)]{AMITA} \\
 \(2\) & \(0\) & \(1\) & \(0\) & \(3^2q\ell\)           & \(3\) & \(0\) & \(1\) & \(1\) & \(1\) & \(2\) & \(4\) & \(\beta\)        & \cite[Thm. 2.5 (5)]{AMITA} \\
 \(2\) & \(0\) & \(2\) & \(0\) & \(3^2\ell_1\ell_2\)    & \(3\) & \(1\) & \(2\) & \(2\) & \(0\) & \(2\) & \(4\) & \(\beta,\gamma\) & Thm. 3.2 (1) \\
 \(1\) & \(1\) & \(0\) & \(0\) & \(3^2\pi\pi^\tau q\)   & \(4\) & \(0\) & \(2\) & \(0\) & \(2\) & \(2\) & \(4\) & \(\alpha,\beta\) & Thm. 3.1 (5) \\
 \(1\) & \(1\) & \(1\) & \(0\) & \(3^2\pi\pi^\tau\ell\) & \(4\) & \(0\) & \(2\) & \(1\) & \(1\) & \(2\) & \(4\) & \(\alpha,\beta\) & Thm. 3.1 (5) \\
 \(1\) & \(1\) & \(0\) & \(1\) & \(3^2\rho\rho^\tau q\) & \(4\) & \(0\) & \(2\) & \(1\) & \(1\) & \(2\) & \(4\) & \(\alpha,\beta\) & Thm. 3.1 (5) \\
 \(3\) & \(0\) & \(0\) & \(0\) & \(3^2q_1q_2q_3\)       & \(4\) & \(0\) & \(2\) & \(0\) & \(3\) & \(3\) & \(8\) & \(\beta\)        & Thm. 3.2 (9) \\
 \(3\) & \(0\) & \(1\) & \(0\) & \(3^2q_1q_2\ell\)      & \(4\) & \(0\) & \(2\) & \(1\) & \(2\) & \(3\) & \(8\) & \(\beta\)        & Thm. 3.2 (9) \\
 \(3\) & \(0\) & \(2\) & \(0\) & \(3^2q\ell_1\ell_2\)   & \(4\) & \(0\) & \(2\) & \(2\) & \(1\) & \(3\) & \(8\) & \(\beta\)        & Thm. 3.2 (9) \\
\hline
\end{tabular}
}
\end{center}
\end{table}

%\newpage
%--------------------------------------------------------------------------------

\renewcommand{\arraystretch}{1.2}

\begin{table}[ht]
\caption{Arithmetic invariants for conductors of Dedekind species 1b}
\label{tbl:Dedekind1b}
\begin{center}
{\normalsize
\begin{tabular}{|c|c|c|c||l||c|c|c||c|c|c|c||c|c|}
\hline
 \(J\) & \(w\) & \(I\) & \(v\) & \(f\)           & \(t\) & \(q^\ast\) & \(r\) & \(G\) & \(B\) & \(T\) & \(m\) & PFT & Reference \\
\hline
 \(1\) & \(0\) & \(0\) & \(0\) & \(3q\)               & \(2\) & \(0\) & \(0\) & \(0\) & \(1\) & \(1\) & \(1\) & \(\beta\)  & \cite[Thm. 2.3 (3)]{AMITA} \\
 \(0\) & \(1\) & \(0\) & \(0\) & \(3\pi\pi^\tau\)     & \(3\) & \(0\) & \(1\) & \(0\) & \(1\) & \(1\) & \(1\) & \(\alpha,\beta\) & \cite[Thm. 2.4 (2)]{AMITA} \\ 
 \(2\) & \(0\) & \(0\) & \(0\) & \(3q_1q_2\)          & \(3\) & \(0\) & \(1\) & \(0\) & \(2\) & \(2\) & \(1\) & \(\beta\)  & \cite[Thm. 2.5 (3)]{AMITA} \\
 \(2\) & \(0\) & \(1\) & \(0\) & \(3q\ell\)           & \(3\) & \(0\) & \(1\) & \(1\) & \(1\) & \(2\) & \(2\) & \(\beta\)  & \cite[Thm. 2.5 (4)]{AMITA} \\
 \(1\) & \(1\) & \(0\) & \(0\) & \(3\pi\pi^\tau q\)   & \(4\) & \(0\) & \(2\) & \(0\) & \(2\) & \(2\) & \(1\) & \(\alpha,\beta\) & Thm. 3.1 (3) \\
 \(1\) & \(1\) & \(1\) & \(0\) & \(3\pi\pi^\tau\ell\) & \(4\) & \(0\) & \(2\) & \(1\) & \(1\) & \(2\) & \(2\) & \(\alpha,\beta\) & Thm. 3.1 (4) \\
 \(1\) & \(1\) & \(0\) & \(1\) & \(3\rho\rho^\tau q\) & \(4\) & \(0\) & \(2\) & \(1\) & \(1\) & \(2\) & \(2\) & \(\alpha,\beta\) & Thm. 3.1 (4) \\
 \(3\) & \(0\) & \(0\) & \(0\) & \(3q_1q_2q_3\)       & \(4\) & \(0\) & \(2\) & \(0\) & \(3\) & \(3\) & \(3\) & \(\beta\)  & Thm. 3.2 (3) \\
 \(3\) & \(0\) & \(1\) & \(0\) & \(3q_1q_2\ell\)      & \(4\) & \(0\) & \(2\) & \(1\) & \(2\) & \(3\) & \(2\) & \(\beta\)  & Thm. 3.2 (4) \\
 \(3\) & \(0\) & \(2\) & \(0\) & \(3q\ell_1\ell_2\)   & \(4\) & \(0\) & \(2\) & \(2\) & \(1\) & \(3\) & \(4\) & \(\beta\)  & Thm. 3.2 (5) \\
\hline
\end{tabular}
}
\end{center}
\end{table}

\newpage
%--------------------------------------------------------------------------------

\renewcommand{\arraystretch}{1.0}

\begin{table}[ht]
\caption{Arithmetic invariants for conductors of Dedekind species 2}
\label{tbl:Dedekind2}
\begin{center}
{\normalsize
\begin{tabular}{|c|c|c|c||l||c|c|c||c|c|c|c||c|c|}
\hline
 \(J\) & \(w\) & \(I\) & \(v\) & \(f\)                         & \(t\) & \(q^\ast\) & \(r\) & \(G\) & \(B\) & \(T\) & \(m\) & PFT & Reference \\
\hline
 \(1\) & \(0\) & \(1\) & \(0\) & \(\ell\)                           & \(1\) & \(1\) & \(0\) & \(1\) & \(0\) & \(1\) & \(1\) & \(\gamma\) & \cite[Thm. 2.3 (2)]{AMITA} \\
 \(0\) & \(1\) & \(0\) & \(1\) & \(\rho\rho^\tau\)                  & \(2\) & \(1\) & \(1\) & \(1\) & \(0\) & \(1\) & \(1\) & \(\alpha,\gamma\) & \cite[Thm. 2.4 (1)]{AMITA} \\
 \(2\) & \(0\) & \(0\) & \(0\) & \(q_1q_2\)                         & \(2\) & \(0\) & \(0\) & \(0\) & \(2\) & \(2\) & \(1\) & \(\beta\)  & \cite[Thm. 2.3 (5)]{AMITA} \\
 \(2\) & \(0\) & \(2\) & \(0\) & \(\ell_1\ell_2\)                   & \(2\) & \(1\) & \(1\) & \(2\) & \(0\) & \(2\) & \(2\) & \(\beta,\gamma\)  & \cite[Thm. 2.5 (2)]{AMITA} \\
 \(1\) & \(1\) & \(0\) & \(0\) & \(\pi\pi^\tau q\)                  & \(3\) & \(0\) & \(1\) & \(0\) & \(2\) & \(2\) & \(1\) & \(\alpha,\beta\) & \cite[Thm. 2.4 (4)]{AMITA} \\
 \(1\) & \(1\) & \(1\) & \(1\) & \(\rho\rho^\tau\ell\)              & \(3\) & \(1\) & \(2\) & \(2\) & \(0\) & \(2\) & \(2\) & \(\alpha,\beta,\gamma\) & Thm. 3.1 (2) \\
 \(3\) & \(0\) & \(0\) & \(0\) & \(q_1q_2q_3\)                      & \(3\) & \(0\) & \(1\) & \(0\) & \(3\) & \(3\) & \(1\) & \(\beta\)  & \cite[Thm. 2.5 (6)]{AMITA} \\
 \(3\) & \(0\) & \(1\) & \(0\) & \(q_1q_2\ell\)                     & \(3\) & \(0\) & \(1\) & \(1\) & \(2\) & \(3\) & \(2\) & \(\beta\)  & \cite[Thm. 2.5 (7)]{AMITA} \\
 \(3\) & \(0\) & \(3\) & \(0\) & \(\ell_1\ell_2\ell_3\)             & \(3\) & \(1\) & \(2\) & \(3\) & \(0\) & \(3\) & \(4\) & \(\beta,\gamma\) & Thm. 3.2 (2) \\
 \(0\) & \(2\) & \(0\) & \(0\) & \(\pi_1\pi_1^\tau\pi_2\pi_2^\tau\) & \(4\) & \(0\) & \(2\) & \(0\) & \(2\) & \(2\) & \(1\) & \(\alpha,\beta\) & Thm. 3.1 (9) \\
 \(2\) & \(1\) & \(0\) & \(0\) & \(\pi\pi^\tau q_1q_2\)             & \(4\) & \(0\) & \(2\) & \(0\) & \(3\) & \(3\) & \(1\) & \(\alpha,\beta\) & Thm. 3.1 (7) \\
 \(2\) & \(1\) & \(1\) & \(0\) & \(\pi\pi^\tau q\ell\)              & \(4\) & \(0\) & \(2\) & \(1\) & \(2\) & \(3\) & \(2\) & \(\alpha,\beta\) & Thm. 3.1 (8) \\ 
 \(2\) & \(1\) & \(0\) & \(1\) & \(\rho\rho^\tau q_1q_2\)           & \(4\) & \(0\) & \(2\) & \(1\) & \(2\) & \(3\) & \(2\) & \(\alpha,\beta\) & Thm. 3.1 (6) \\
 \(4\) & \(0\) & \(0\) & \(0\) & \(q_1q_2q_3q_4\)                   & \(4\) & \(0\) & \(2\) & \(0\) & \(4\) & \(4\) & \(3\) & \(\beta\)  & Thm. 3.2 (6) \\
 \(4\) & \(0\) & \(1\) & \(0\) & \(q_1q_2q_3\ell\)                  & \(4\) & \(0\) & \(2\) & \(1\) & \(3\) & \(4\) & \(2\) & \(\beta\)  & Thm. 3.2 (7) \\
 \(4\) & \(0\) & \(2\) & \(0\) & \(q_1q_2\ell_1\ell_2\)             & \(4\) & \(0\) & \(2\) & \(2\) & \(2\) & \(4\) & \(4\) & \(\beta\)  & Thm. 3.2 (8) \\
\hline
\end{tabular}
}
\end{center}
\end{table}

%\newpage
%--------------------------------------------------------------------------------

\noindent
The simultaneous proof of
\cite[Thm., \S\ 1, p. 8]{Ho1971},
\cite[Thm. 1.1, p. 251]{AMITA},
and Theorem
\ref{thm:Rank2}
in the Tables
\ref{tbl:Dedekind1a},
\ref{tbl:Dedekind1b} and
\ref{tbl:Dedekind2},
which are already purged from \textit{nilets}, i.e. multiplets with \(m=0\), namely
\(f\in\lbrace 3\ell, 3\rho\rho^\tau, 3\ell_1\ell_2, 3\rho\rho^\tau\ell, 3\ell_1\ell_2\ell_3\rbrace\) of species 1b,
and \(f\in\lbrace q, \pi\pi^\tau, q\ell, \pi\pi^\tau\ell,\rho\rho^\tau q, q\ell_1\ell_2\rbrace\) of species 2,
establishes the rank \(r\) of the ambiguous \(3\)-class group \(Cl_3^{(\sigma)}(N)\)
which is only an approximation of the entire \(3\)-class group \(Cl_3(N)\).

%\noindent
The systematic investigation of the full \(3\)-class group \(Cl_3(N)\)
of the normal closure \(N\) of arbitrary non-Galois cubic number fields \(L\)
began in \(1933\) with conclusions from the class number relation
\(h_N=\frac{u}{3^\varrho}\cdot h_L^2\cdot h_K\)
of Arnold Scholz
\cite[p. 213, p. 216]{So1933},
involving the quadratic subfield \(K\) of \(N\),
the torsion free Dirichlet unit rank \(\varrho\) of \(L\),
and the index \(u=(E_N:E_0)\) of the subgroup \(E_0\) generated by the units of all subfields of \(N\)
in the unit group \(E_N\) of \(N\).
For the situation of a pure cubic field \(L=\mathbb{Q}(\sqrt[3]{D})\) with \(\varrho=1\) and \(h_K=1\),
the index \(u\) can only take two values,
\(u=1\) for the \textit{differential principal factorization type} (briefly PFT) \(\alpha\) in the sense of
\cite[Thm. 2.1, p. 254]{AMITA},
and \(u=3\) for the types \(\beta\) and \(\gamma\).
In particular, the isomorphism \(Cl_3(N)\simeq Cl_3(L)\times Cl_3(L^\sigma)\) to the \textit{direct product}
is only possible for the types \(\beta\) and \(\gamma\)
\cite[p. 219]{So1933},
since, due to the relative principal factorization
\cite{BC1971}
in a field \(N\) of type \(\alpha\), the direct product shrinks to a quotient
\((Cl_3(L)\times Cl_3(L^\sigma))/(\mathbb{Z}/3\mathbb{Z})\).
This was confirmed later by Gerth
\cite[Thm. 5.1, p. 61]{Ge1975}
in the important case \(w=0\), that is,
if the conductor \(f\) is only divisible by non-split primes \(q_i\), \(1\le i\le J\).
 
%\noindent
The need for a detailed description of the \(3\)-class group \(Cl_3(N)\)
of the normal closure \(N\) of a pure cubic field \(L\),
in terms of power products of prime ideals whose classes generate the group,
has its origin in the \(1992\) Doctoral Thesis of Ismaili
\cite{Is},
where he studied the capitulation kernels \(\ker(T_i)\)
in an elementary bicyclic \(3\)-class group \(Cl_3(N)\)
with respect to the transfers \(T_i:\,Cl_3(N)\to Cl_3(K_i)\) from \(N\)
to the four unramified cyclic cubic extensions \(K_1,\ldots,K_4\) of \(N\).
Together with El Mesaoudi
\cite{IsEM},
he determined necessary and sufficient conditions for \(N\)
to possess an elementary bi-homo-cyclic \(3\)-class group
\(Cl_3(N)\cong\mathbb{Z}/3\mathbb{Z}\times\mathbb{Z}/3\mathbb{Z}\),
isomorphic to the direct product of two copies of \(Cl_3(L)\cong\mathbb{Z}/3\mathbb{Z}\),
which enforces \(u=3\), i.e. PFT \(\beta\) or \(\gamma\).

In
\cite{ATMI},
we studied the capitulation
in a non-elementary bi-hetero-cyclic \(3\)-class group
\(Cl_3(N)\cong\mathbb{Z}/9\mathbb{Z}\times\mathbb{Z}/3\mathbb{Z}\),
based on a general theory of such class groups in
\cite{Ma2023}.
This scenario occurs for the conductor \(f=\rho\rho^\tau\),
i.e. a prime \(p\equiv 1\,(\mathrm{mod}\,9)\)
\cite[Thm. 2.4 (1)]{AMITA},
with \(u=1\) and \(Cl_3(L)\cong\mathbb{Z}/9\mathbb{Z}\).

In a forthcoming work
\cite{Ao},
we shall investigate the capitulation
in an elementary tri-homo-cyclic \(3\)-class group
\(Cl_3(N)\cong\mathbb{Z}/3\mathbb{Z}\times\mathbb{Z}/3\mathbb{Z}\times\mathbb{Z}/3\mathbb{Z}\),
which enforces \(u=1\), i.e. PFT \(\alpha\),
for parity reasons in the class number formula of Scholz.

\newpage
%-----------------------------------------------------------------------

\section{Conductors and their multiplicities}
\label{s:Context}

Before we give the proof of Theorem \ref{thm:Rank2} in \S\ \ref{s:Proof},
we devote the present section to an overview of the special properties
of pure cubic fields $L=\mathbb{Q}(\sqrt[3]{D})$
possessing the radicands $D$ of the shapes listed in equation \eqref{eqn:Rank2} of the introduction \S\ \ref{s:Intro}.

%-----------------------------------------------------------------------

\subsection{Conductors with splitting prime divisors}
\label{ss:Split}

Let us start with giving more details concerning the leading eight lines
of equation \eqref{eqn:Rank2} in our Theorem \ref{thm:Rank2},
where $D$ is divisible by a prime $p_1\equiv 1\,(\mathrm{mod}\,3)$
which splits in $K$.

\begin{theorem}
\label{thm:Ismaili1}
Let the conductor of $N/K$ be $f=3^\varepsilon\cdot p_1\cdots p_w\cdot q_1\cdots q_J$
as in Equation
\eqref{eqn:PrimeDecCond}
with $0\le\varepsilon\le 2$, $T=w+J\ge 1$,
and pairwise distinct primes
$p_{i}\equiv 1\,(\mathrm{mod}\,3)$ for $1\le i\le w$,
and $q_{i}\equiv -1\,(\mathrm{mod}\,3)$ for $1\le i\le J$.
Briefly denote the multiplicity of $f$ by $m:=m(f)$.
Assume that $w\ge 1$.
Then, $\operatorname{rank}\,(Cl_{3}^{(\sigma)}(N))=2$
$\Longleftrightarrow$ $L$ belongs to one of the following multiplets:

\begin{enumerate}
\item[$(1)$]
doublets, $m=2$, of type $(\alpha^x,\beta^y,\gamma^z)$, $x+y+z=2$,
such that $f=9p_1$ with $p_1\equiv 1\,(\mathrm{mod}\,9)$,
\item[$(2)$]
doublets, $m=2$, of type $(\alpha^x,\beta^y,\gamma^z)$, $x+y+z=2$,
such that $f=p_1q_1$ with $p_1\equiv 1\,(\mathrm{mod}\,9)$ and $q_1\equiv 8\,(\mathrm{mod}\,9)$,
\item[$(3)$]
singlets, $m=1$, of type $\alpha$ or $\beta$
such that $f=3p_1q_1$ with $p_1\equiv 4,7\,(\mathrm{mod}\,9)$, $q_1\equiv 2,5\,(\mathrm{mod}\,9)$,
\item[$(4)$]
doublets, $m=2$, of type $(\alpha^x,\beta^y)$, $x+y=2$, such that $f=3p_1q_1$ with \\
either $(p_1\equiv 1\,(\mathrm{mod}\,9)$ and $q_1\equiv 2,5\,(\mathrm{mod}\,9))$
or $(p_1\equiv 4,7\,(\mathrm{mod}\,9)$ and $q_1\equiv 8\,(\mathrm{mod}\,9))$,
\item[$(5)$]
quartets, $m=4$, of type $(\alpha^x,\beta^y)$, $x+y=4$, 
such that $f=9p_1q_1$ with $p_1\not\equiv 1\,(\mathrm{mod}\,9)$ or $q_1\not\equiv 8\,(\mathrm{mod}\,9)$,
\item[$(6)$]
doublets, $m=2$, of type $(\alpha^x,\beta^y)$, $x+y=2$, 
such that $f=p_1q_1q_2$ with $p_1\equiv 1\,(\mathrm{mod}\,9)$ and $q_1,q_2\equiv 2,5\,(\mathrm{mod}\,9)$,
\item[$(7)$]
singlets, $m=1$, of type $\alpha$ or $\beta$
such that $f=p_1q_1q_2$ with $p_1\equiv 4,7\,(\mathrm{mod}\,9)$ and $q_1,q_2\equiv 2,5\,(\mathrm{mod}\,9)$,
\item[$(8)$]
doublets, $m=2$, of type $(\alpha^x,\beta^y)$, $x+y=2$,
such that $f=p_1q_1q_2$ with $p_1\equiv 4,7\,(\mathrm{mod}\,9)$, $q_1\equiv 2,5\,(\mathrm{mod}\,9)$ and $q_2\equiv 8\,(\mathrm{mod}\,9)$,
\item[$(9)$]
singlets, $m=1$, of type $\alpha$ or $\beta$
such that $f=p_1p_2$ with $p_1,p_2\equiv 4,7\,(\mathrm{mod}\,9)$.
\end{enumerate}

\noindent
There exist infinitely many multiplets with conductors of all these shapes $(1)$--$(9)$.
\end{theorem}

%--------------------------------------------------------------------------------

\begin{proof}
The multiplicity $m(f)$ of each conductor is calculated by means of Formula
\eqref{eqn:Multiplicity},
using the sequence $(X_k)_{k\ge -1}=(\frac{1}{2},0,1,1,3,5,11,\ldots)$:
\begin{enumerate}
\item[$(1)$]
For  $D=3^{e}p_1^{e_1} \not\equiv\pm 1\,(\mathrm{mod}\,9)$
with $p_1 \equiv 1\,(\mathrm{mod}\,9)$,
we have $f=3^2p_1$ of species 1a.
We must take into consideration that $G=1$, $B=0$, $T=1$,
where $G$, $B$ and $T$ are the numbers defined in \eqref{eqn:Multiplicity},
and we obtain $m(f)=2^T=2$, a doublet, independently of $G$ and $B$. \\
 Since
$\mathbb{Q}(\sqrt[3]{ab^2})=\mathbb{Q}(\sqrt[3]{a^2b})$ for $  a,
b \in \mathbb{N}$, we can choose $e_1=1$, so $D=3^{e}p$ with $e\in
\lbrace 1,2 \rbrace$. As  $ p \equiv 1  \pmod  3$, then
  $p=\pi_1 \pi_2$, where $\pi_1$ and $\pi_2$ are primes of $K$
  such that $\pi_1^{\tau}=\pi_2$ and $\pi_1 \equiv \pi_2 \equiv 1 ~(\bmod~3\mathcal{O}_{K})$, and $p$ is totally ramified in $L$. So  $\pi_1\mathcal{O}_{N}=\mathcal{P}_1^3$
and $\pi_2 \mathcal{O}_{N}=\mathcal{P}_2^3$ where $\mathcal{P}_1, \mathcal{P}_2$
are two prime ideals of $N$. The fact that $D \not \equiv \pm 1 ~(\bmod~9) $ implies that $3$ is totally ramified in $L$, then $\lambda$ is ramified in $N/K$.
  So the number of ideals which are ramified in $N/K$ is $t=3.$ As $3=-\zeta_3^2\lambda^2$, then
$N=K(\sqrt[3]{x})$, where $x=\zeta_3^2\lambda^2
\pi_1 \pi_2$. If $
p \equiv 1 \pmod 9$, then $\pi_1$ and $\pi_2$ are congruent to
$ 1 ~(\bmod~\lambda^3)$, so by \cite[Lem. 3.3, p. 264]{AMITA},  $\zeta_3$ is a norm of an element of $N-\lbrace 0\rbrace$ and $q^{*}=1$,
  thus $\operatorname{rank}\,(Cl_{3}^{(\sigma)}(N))=2$.
  The opposite implication of the rank condition $\operatorname{rank}\,(Cl_{3}^{(\sigma)}(N))=2$
to this shape of conductor will be proved in \S\S\ \ref{01}.
%%%%%%%%
\item[$(2)$]
For $D=p_1^{e_1}q_1^{f_1}$ with $p_1\equiv -q_1\equiv 1\,(\mathrm{mod}\,9)$,
we have $f=p_1q_1$ of species 2.
We must take into consideration that $G=2$, $B=0$,
and we obtain $m(f)=2^G\cdot X_{B-1}=2^2\cdot\frac{1}{2}=2$, a doublet. \\
Here $e_1,f_1 \in\{1,2\}$. Since
$\mathbb{Q}(\sqrt[3]{ab^2})=\mathbb{Q}(\sqrt[3]{a^2b})$ for $  a,
b \in \mathbb{N}$, we can choose $e_1=1$, so $D=pq^{f_1} $. As $D \equiv \pm1 \pmod  9$,
 $3$ is not ramified in the field $L$, and then $\lambda$ is not ramified in $N/K$.
 the fact that  $ p \equiv 1  \pmod  3$, implies that $p=\pi_1 \pi_2$, where $\pi_1$ and $\pi_2$ are primes of $K$
 such that  $\pi_1^{\tau}=\pi_2$ and $\pi_1 \equiv \pi_2 \equiv 1 ~(\bmod~3\mathcal{O}_{K})$,
 since $p$ is totally ramified in $L$, then $\pi_1$ and $\pi_2$ are totally ramified in $N$.
 As  $ q \equiv -1  \pmod  3$, $q$ is inert in $K$. Then, the primes ramified in $N/K$
 are $\pi_1, \pi_2$ and $q$, and $t=3$. Since $ p \equiv -q \equiv 1  \pmod 9$, then $\pi_1 \equiv \pi_2 \equiv \pi \equiv 1 ~(\bmod~\lambda^3)$,
  where $-q=\pi$ is a prime number of $K$. Then,
 $N=K(\sqrt[3]{x})$, where $x=\pi_1^{e_1} \pi_2^{e_1} \pi^{f_1}$.  The fact that all primes $\pi$, $\pi_1$ and $\pi_2$ are congruent to $ 1 ~(\bmod~\lambda^3)$,
  implies by  \cite[Lem. 3.3, p. 264]{AMITA} that $\zeta_3$ is a norm of an element of $N-\lbrace 0\rbrace$ and $q^{*}=1$. Hence, $\operatorname{rank}\,(Cl_{3}^{(\sigma)}(N))=2$.
 The opposite implication of the rank condition $\operatorname{rank}\,(Cl_{3}^{(\sigma)}(N))=2$
to this shape of conductor will be proved in \S\ \ref{03}.
%%%%%%%%
\item[$(3)$]
For $D=p_1^{e_1}q_1^{f_1}\not\equiv\pm 1\,(\mathrm{mod}\,9)$
with $(p_1\equiv 4,7\,(\mathrm{mod}\,9)$ and $q_1\equiv 2,5\,(\mathrm{mod}\,9))$,
we have $f=3p_1q_1$ of species 1b.
We must take into consideration that $G=0$, $B=2$,
and we obtain $m(f)=2^G\cdot X_{B}=1\cdot 1=1$, a singlet. \\
Here $e_1,f_1 \in\{1,2\}$. As
$\mathbb{Q}(\sqrt[3]{ab^2})=\mathbb{Q}(\sqrt[3]{a^2b})$ $ \forall a,
b \in \mathbb{N}$, we can choose $e_1=1$, i.e $D=pq^{f_1} $. Since $D \not \equiv \pm1 \pmod  9$,
 then $3$ is ramified in the field $L$, so $\lambda=1-\zeta_3$ is  ramified in $N/K$.
 Since  $ p \equiv 1  \pmod  3$, then $p=\pi_1 \pi_2$, where $\pi_1$ and $\pi_2$ are primes of $K$
 such that  $\pi_1^{\tau}=\pi_2$ and $\pi_1 \equiv \pi_2 \equiv 1 ~(\bmod~3\mathcal{O}_{K})$,
 the prime $p$ is totally ramified in $L$, so $\pi_1$ and $\pi_2$ are totally ramified in $N$.
 Since  $ q \equiv -1  \pmod  3$, then $q$ is inert in $K$. Thus the primes ramified in $N/K$
 are $\pi_1, \pi_2$ and $q$, and $t=4$.\\
  Let $-q=\pi$ be a prime number of $K$, and put $x=\pi_1^{e_1} \pi_2^{e_1} \pi^{f_1}$, then
 $N=K(\sqrt[3]{x})$. Then:
 \begin{itemize}
 \item If $ p \not \equiv 1  \pmod  9,$  then by  \cite[Lem. 3.1, p. 263]{AMITA},  $\pi_1 \not \equiv  1 ~(\bmod~\lambda^3)$ and  $ \pi_2 \not \equiv 1 ~(\bmod~\lambda^3)$, then according to  \cite[Lem. 3.3, p. 264]{AMITA}, $\zeta_3$ is not a norm of an element of $N-\lbrace 0\rbrace$.
 \item If  $ q \not \equiv -1 (\bmod
   9),$  then by  \cite[Lem. 3.2, p. 264]{AMITA},  $\pi \not \equiv  1 ~(\bmod~\lambda^3)$, then according to  \cite[Lem. 3.3, p. 264]{AMITA}, $\zeta_3$ is not a norm of an element of $N-\lbrace 0\rbrace$.
 \end{itemize}
In all cases we have $q^{*}=0$.
  We conclude that $\operatorname{rank}\,(Cl_{3}^{(\sigma)}(N))=2$.
 The opposite implication of the rank condition $\operatorname{rank}\,(Cl_{3}^{(\sigma)}(N))=2$
to this shape of conductor will be proved in \S\ \ref{03}.
%%%%%%%%
\item[$(4)$]
For $D=p_1^{e_1}q_1^{f_1}\not\equiv\pm 1\,(\mathrm{mod}\,9)$
with either $(p_1\equiv 1\,(\mathrm{mod}\,9)$ and $q_1\equiv 2,5\,(\mathrm{mod}\,9))$
or $(p_1\equiv 4,7\,(\mathrm{mod}\,9)$ and $q_1\equiv 8\,(\mathrm{mod}\,9))$,
we have $f=3p_1q_1$ of species 1b.
We must take into consideration that $G=1$, $B=1$,
and we obtain $m(f)=2^G\cdot X_{B}=2\cdot 1=2$, a doublet.
%%%%%%%%
\item[$(5)$]
For $D=3^{e}p_1^{e_1}q_1^{f_1}\not\equiv\pm 1\,(\mathrm{mod}\,9)$
with $p_1\not\equiv 1\,(\mathrm{mod}\,9) \text{ or } q_1\not\equiv 8\,(\mathrm{mod}\,9)$,
we have $f=3^2p_1q_1$ of species 1a, $T=2$,
and we get $m(f)=2^T=2^2=4$, a quartet, independently of $G$ and $B$. \\
Here $e_1,f_1\in\{1,2\}$. As
$\mathbb{Q}(\sqrt[3]{ab^2})=\mathbb{Q}(\sqrt[3]{a^2b})$ $ \forall a,
b \in \mathbb{N}$, we can choose $e_1=1$, i.e $D=3^{e}pq^{f_1}$. Since $D \not \equiv \pm1 \pmod  9$, we reason as in case (3) and we obtain $t=4$. Since $3=-\zeta_3^2\lambda^2$,  then
 $N=K(\sqrt[3]{x})$, where $x=\zeta_3^2\lambda^2\pi_1^{e_1} \pi_2^{e_1} \pi^{f_1}$. As in case (3), if  $ p \not \equiv 1  \pmod  9 \ or \ q \not \equiv -1 \pmod 9, $ then by using  \cite[Lem. 3.1, Lem. 3.2, Lem. 3.3, p. 263-264]{AMITA}, we get $\zeta_3$ is not a norm of an element of $N-\lbrace 0\rbrace$ and $q^{*}=0$.
  We conclude that $\operatorname{rank}\,(Cl_{3}^{(\sigma)}(N))=2$. The opposite implication of the rank condition $\operatorname{rank}\,(Cl_{3}^{(\sigma)}(N))=2$
to this shape of conductor will be proved in \S\ \ref{03}.
%%%%%%%%
\item[$(6)$]
For $D=p_1^{e_1}q_1^{f_1}q_2^{f_2}\equiv\pm 1\,(\mathrm{mod}\,9)$
with $p_1\equiv 1\,(\mathrm{mod}\,9)$ and $q_1,q_2\equiv 2 \text{ or } 5\,(\mathrm{mod}\,9)$, 
we have $f=p_1q_1q_2$ of species 2.
We must take into consideration that $G=1$, $B=2$,
and we obtain $m(f)=2^G\cdot X_{B-1}=2\cdot 1=2$, a doublet. \\
Here $e_1,f_1 \in\{1,2\}$. As
$\mathbb{Q}(\sqrt[3]{ab^2})=\mathbb{Q}(\sqrt[3]{a^2b})$ $ \forall a,
b \in \mathbb{N}$, we can choose $e_1=1$, i.e $D=pq_1^{f_1}q_2^{f_2} $. Since $D \equiv \pm1 \pmod  9$,
 then $3$ is not ramified in the field $L$, so $\lambda$ is not ramified in $N/K$.
 Since  $ p \equiv 1  \pmod  3$, then $p=\pi_1 \pi_2$, where $\pi_1$ and $\pi_2$ are primes of $K$
 such that  $\pi_1^{\tau}=\pi_2$ and $\pi_1 \equiv \pi_2 \equiv 1 ~(\bmod~3\mathcal{O}_{K})$,
 the prime $p$ is totally ramified in $L$, so $\pi_1$ and $\pi_2$ are totally ramified in $N$.
 Since  $ q_1 \equiv q_2 \equiv -1  \pmod  3$, then $q_1$ and $q_2$ are inert in $K$. Thus the primes ramified in $N/K$
 are $\pi_1, \pi_2$, $q_1$, and $q_2$. Then, $t=4$.\\
  Let $-q_1=\pi'$ and $-q_2=\pi''$ be two prime numbers of $K$, and put $x=\pi_1^{e_1} \pi_2^{e_1} \pi^{'f_1}\pi^{''f_2}$, then
 $N=K(\sqrt[3]{x})$. As $q_1 \not \equiv -1 (\bmod 9),$ then by  \cite[Lem. 3.2, p. 264]{AMITA},  $\pi' \not \equiv  1 ~(\bmod~\lambda^3)$, then according to  \cite[Lem. 3.3, p. 264]{AMITA}, $\zeta_3$ is not a norm of an element of $N-\lbrace 0\rbrace$, so $q^{*}=0$. Thus, $\operatorname{rank}\,(Cl_{3}^{(\sigma)}(N))=2$.
 The opposite implication of the rank condition $\operatorname{rank}\,(Cl_{3}^{(\sigma)}(N))=2$
to this shape of conductor will be proved in \S\ \ref{05}.
%%%%%%%%
\item[$(7)$]
For $D=p_1^{e_1}q_1^{f_1}q_2^{f_2}\equiv\pm 1\,(\mathrm{mod}\,9)$
with $p_1,-q_1,-q_2\equiv 4 \text{ or } 7\,(\mathrm{mod}\,9)$,
we have $f=p_1q_1q_2$ of species 2.
We must take into consideration that $G=0$, $B=3$,
and we obtain $m(f)=2^G\cdot X_{B-1}=1\cdot 1=1$, a singlet. \\
We reason as in case (6), we get $t=4$ and $q^{*}=0$, then $\operatorname{rank}\,(Cl_{3}^{(\sigma)}(N))=2$. 
 The opposite implication of the rank condition $\operatorname{rank}\,(Cl_{3}^{(\sigma)}(N))=2$
to this shape of conductor will be proved in \S\ \ref{05}.
%%%%%%%%
\item[$(8)$]
For $D=p_1^{e_1}q_1^{f_1}q_2^{f_2}\equiv\pm 1\,(\mathrm{mod}\,9)$
with $p_1,-q_2\equiv 4 \text{ or } 7\,(\mathrm{mod}\,9)$
and $q_1\equiv 8\,(\mathrm{mod}\,9)$,
we have $f=p_1q_1q_2$ of species 2.
We must take into consideration that $G=1$, $B=2$,
and we obtain $m(f)=2^G\cdot X_{B-1}=2\cdot 1=2$, a doublet. \\
 We reason as in case (6), we get $t=4$ and $q^{*}=0$, then $\operatorname{rank}\,(Cl_{3}^{(\sigma)}(N))=2$.
 The opposite implication of the rank condition $\operatorname{rank}\,(Cl_{3}^{(\sigma)}(N))=2$
to this shape of conductor will be proved in \S\ \ref{05}.
%%%%%%%%
\item[$(9)$]
For $D=p_1^{e_1}p_2^{e_2}\equiv\pm 1\,(\mathrm{mod}\,9)$
such that $p_{i}\not\equiv 1\,(\mathrm{mod}\,9)$ for some $i\in\lbrace 1,2\rbrace$,
we have $f=p_1p_2$ of species 2.
We must take into consideration that $G=0$, $B=2$,
and we obtain $m(f)=2^G\cdot X_{B-1}=1\cdot 1=1$, a singlet. \\
Here $e_1, e_2 \in   \{1,2\}$. We have  $p_1=\pi_1\pi_2$ and $p_2=\pi_3\pi_4$, where $\pi_1$, $\pi_2$, $\pi_3$ and
$\pi_4$ are primes of $K$ such that
$\pi_2=\pi_1^{\tau}$, $\pi_4=\pi_3^{\tau}$, and $\pi_1 \equiv \pi_2 \equiv \pi_3 \equiv \pi_4 \equiv 1 ~(\bmod~3\mathcal{O}_{K})$,
and since $p_1$ and $p_2$ are ramified in $L$, then $\pi_1$, $ \pi_2$, $ \pi_3$, and $ \pi_4$ are ramified in $N$. The fact that $D  \equiv \pm1 \pmod  9$ implies that $3$ is  not ramified in $L$, so $\lambda$ is not  ramified in $N/K$. Then we get $t=4$. \\
As $ p_1 \not \equiv 1 \pmod 9 $ or $ p_2 \not \equiv 1 \pmod 9 $, then by  \cite[Lem. 3.1, p. 263]{AMITA}, $ \exists i \in \lbrace 1, 2, 3, 4 \rbrace$ such that $\pi_{i}$ is not congruent to $ 1 ~(\bmod~\lambda^3)$, then according to \cite[Lem. 3.3, p. 264]{AMITA}, $\zeta_3$ is not norm of an element of $N-\lbrace 0\rbrace$ and $q^{*}=0$. Hence, $\operatorname{rank}\,(Cl_{3}^{(\sigma)}(N))=2.$
 The opposite implication of the rank condition $\operatorname{rank}\,(Cl_{3}^{(\sigma)}(N))=2$
to this shape of conductor will be proved in \S\ \ref{06}.
\end{enumerate}

\noindent
The principal factorization type
is a consequence of the estimates in \cite[Thm. 2.1, p. 254]{AMITA}.
Since $s=1$, we have $0\le R\le s=1$ and type $\alpha$ with relative PF may generally occur.
For all items $(1)$--$(9)$, we have $3\le t\le 4$
and thus $1\le A\le\min(2,t-s)=2$ with $2=3-1\le t-s\le 4-1=3$,
which generally enables type $\beta$ with absolute PF,
but type $\gamma$ may occur.
For items $(3)$--$(9)$, there exists either a prime factor $p_1\equiv 4,7\,(\mathrm{mod}\,9)$
or a prime factor $q_1\equiv 2,5\,(\mathrm{mod}\,9)$,
whence type $\gamma$ is impossible,
because $\zeta_3$ can be norm of a unit in $N$ only if the prime factors of $f$
are $3$ or $\ell_j\equiv 1,8\,(\mathrm{mod}\,9)$.

All claims on the infinitude of the various sets of conductors $f$
are a consequence of Dirichlet's theorem on primes in arithmetic progressions.
\end{proof}

%\newpage
%-----------------------------------------------------------------------

\noindent
\textbf{Examples 1.} In all of our applications,
we present the structure of the $3$-class groups of pure cubic fields $L$ and of their Galois closures $N$
in the case where the conductor $f$ contains splitting prime divisors.
Computations were performed with Magma
\cite{BCP1997,BCFS2023,MAGMA2023}
and Pari/GP
\cite{PARI2023}.
 
Let $Cl_{3}(N)$ (respectively $Cl_{3}(L)$) be the \(3\)-class group  of \(N\) (respectively of $L$),
and $N_3^{(1)}$ be the maximal abelian unramified $3$-extension of $N$. 
Let $u$ be the index of the subgroup generated by the units of intermediate fields of the extension $N/\mathbb{Q}$ in the group of units of $N$.
According to \cite[\S\ 12, Theorem 12.1, p. 229]{BC1971}, there are two possibilities, either $u=1$ or $u=3$.

%\paragraph{}
%%%%%%%%%%%%%%%%%%%%%%%%%%%%%%%%%%%%%%%%%%%%%%%%%

In Table \ref{tb1}, we start with conductors $f=9p_1$ of species 1a, where $p_1\equiv 1\,(\mathrm{mod}\,9)$.
See Theorem
\ref{thm:Ismaili1},
item (1), and the first line of Equation
\eqref{eqn:Rank2}.

\renewcommand{\arraystretch}{1.0}
 
\begin{table}[ht]
\caption{Radicands $D=3^ep_1$ of species 1a with $m=2$}
\label{tb1}
\begin{center}
\begin{tabular}{|rlrr|ccc|cc|}
\hline
 $D$    & $3^e$ & $p_1$ & $f$    & $\left(\dfrac{3}{p_1}\right)_3$ & $u$ & PFT & $Cl_{3}(L)$ & $Cl_{3}(N)$ \\
\hline
   $57$ & $3$   &  $19$ &  $171$ & $\neq 1$ & $3$ & $\gamma$ & $(3)$    &   $(3,3)$   \\
  $111$ & $3$   &  $37$ &  $333$ & $\neq 1$ & $3$ & $\gamma$ & $(3)$    &   $(3,3)$   \\
  $219$ & $3$   &  $73$ &  $657$ & $=1$     & $3$ & $\beta$  & $(3,3)$  & $(9,3,3)$   \\
  $327$ & $3$   & $109$ &  $981$ & $\neq 1$ & $3$ & $\gamma$ & $(3)$    &   $(3,3)$   \\
  $813$ & $3$   & $271$ & $2439$ & $=1$     & $1$ & $\alpha$ & $(27,3)$ & $(27,27,3)$ \\
  $921$ & $3$   & $307$ & $2763$ & $=1$     & $3$ & $\beta$  & $(3,3)$  & $(9,3,3)$   \\
 $1569$ & $3$   & $523$ & $4707$ & $=1$     & $1$ & $\alpha$ & $(9,3)$  & $(9,9,3)$   \\
 $1629$ & $3^2$ & $181$ & $1629$ & $\neq 1$ & $3$ & $\gamma$ & $(3)$    &   $(3,3)$   \\
 $1791$ & $3^2$ & $199$ & $1791$ & $\neq 1$ & $3$ & $\gamma$ & $(3)$    &   $(3,3)$   \\
\hline
\end{tabular}
\end{center}
\end{table}

\begin{remark}
\label{rem1}
In Table \ref{tb1}, if the $3$-class group $Cl_{3}(N)$ is of type $(3,3)$, then $N^{\ast}=N_3^{(1)}$.
In this case, the cubic residue symbol $\left(\dfrac{3}{p_1}\right)_3$ is different from $1$
if and only if $3$ divides exactly the class number of $L$, which implies that $u=3$.
Consequently, Table \ref{tb1} confirms the previous results by Ismaili and El Mesaoudi
\cite[Theorem 1, case $(i)$, p. 157; Theorem 2 and Corollary 1, pp. 161--162]{IsEM}.
\end{remark}

%\newpage
%\paragraph{}
%%%%%%%%%%%%%%%%%%%%%%%%%%%%%%%%%%%%%%%%%%%%%%%%%%%%%

\noindent
Table \ref{tb2} deals with conductors $f=p_1q_1$ of species 2, where $p_1\equiv -q_1\equiv 1\,(\mathrm{mod}\,9)$.
See Theorem
\ref{thm:Ismaili1},
item (2), and the second line of Equation
\eqref{eqn:Rank2}.

\renewcommand{\arraystretch}{1.0}
 
\begin{table}[ht]
\caption{Radicands $D=p_1q_1$ of species 2 with $m=2$}
\label{tb2}
\begin{center}
\begin{tabular}{|rrrr|cccc|cc|}
\hline
 $D$    & $p_1$ & $q_1$ & $f$    & $\left(\dfrac{3}{p_1}\right)_3$ & $\left(\dfrac{q_{1}}{p_1}\right)_3$ & $u$ & PFT & $Cl_{3}(L)$ & $Cl_{3}(N)$ \\
\hline
  $323$ &  $19$ &  $17$ &  $323$ & $\neq 1$ & $\neq 1$ & $3$ & $\gamma$ & $(3)$    & $(3,3)$      \\
 $1007$ &  $19$ &  $53$ & $1007$ & $\neq 1$ & $\neq 1$ & $3$ & $\gamma$ & $(3)$    & $(3,3)$      \\
 $1241$ &  $73$ &  $17$ & $1241$ & $=1$     & $=1$     & $1$ & $\alpha$ & $(9,3)$  & $(9,3,3,3)$  \\
 $1349$ &  $19$ &  $71$ & $1349$ & $\neq 1$ & $\neq 1$ & $3$ & $\gamma$ & $(3)$    & $(3,3)$      \\
 $2033$ &  $19$ & $107$ & $2033$ & $\neq 1$ & $=1$     & $3$ & $\beta$  & $(3,3)$  & $(9,3,3)$    \\ 
 $3401$ &  $19$ & $179$ & $3401$ & $\neq 1$ & $=1$     & $1$ & $\alpha$ & $(27,3)$ & $(27,9,3,3)$ \\
 $3869$ &  $73$ &  $53$ & $3869$ & $=1$     & $\neq 1$ & $3$ & $\gamma$ & $(3)$    & $(3,3)$      \\
\hline
\end{tabular}
\end{center}
\end{table}

\begin{remark}
\label{rem2}
In Table \ref{tb2}, if the $3$-class group $Cl_{3}(N)$ is of type $(3,3)$, then $N^{\ast}= N_3^{(1)}$.
In this case, the structure of the $3$-class group $Cl_{3}(N)$ is
independent of the cubic residue symbol $\left(\dfrac{3}{p_1}\right)_3$.
However, the cubic residue symbol $\left(\dfrac{q_{1}}{p_1}\right)_3$ is different from $1$
if and only if $3$ divides exactly the class number of $L$, which implies that $u=3$.
Consequently, Table \ref{tb2} confirms the previous results by Ismaili and El Mesaoudi
\cite[Theorem 1, case $(ii)$, p. 157; Theorem 3 and Corollary 2, pp. 162--163]{IsEM}.
\end{remark}

%\paragraph{}
%%%%%%%%%%%%%%%%%%%%%%%%%%%%%%%%%%%%%%%%%%%%%%%%%%%%%

Table \ref{tb3} concerns conductors $f=3p_1q_1$ of species 1b, where $p_1\equiv 1\,(\mathrm{mod}\,9)$ and $q_1\equiv 2,5\,(\mathrm{mod}\,9)$.
See Theorem
\ref{thm:Ismaili1},
item (4), and the third line of Equation
\eqref{eqn:Rank2}.
Other examples are presented in
\cite[Exm. 7.4, pp. 116--117]{AMI}:
\(14801=19^2\cdot 41\) and \(56129=37^2\cdot 41\) contain $p_1\equiv 1\,(\mathrm{mod}\,9)$,
whereas \(833=7^2\cdot 17\), \(8959=17^2\cdot 31\), \(97997=43^2\cdot 53\) contain $q_1\equiv 8\,(\mathrm{mod}\,9)$,
also with \(m=2\).
However, \(1573=11^2\cdot 13\), \(4901=13^2\cdot 29\), \(22747=23^2\cdot 43\), \(32269=23^2\cdot 61\) are singlets \((m=1)\)
belonging to Theorem
\ref{thm:Ismaili1},
item (3).
All of them are M0-fields.
(Cfr.
\cite[Thm. 4--6]{MaSl2023}.)

\renewcommand{\arraystretch}{1.0}
 
\begin{table}[ht]
\caption{Radicands $D=p_1q_1$ of species 1b with $m=2$}
\label{tb3}
\begin{center}
\begin{tabular}{|rrrr|ccc|cc|}
\hline
 $D$    & $p_1$ & $q_1$ & $f$    & $\left(\dfrac{q_{1}}{p_1}\right)_3$ & $u$ & PFT & $Cl_{3}(L)$ & $Cl_{3}(N)$    \\
\hline
   $38$ &  $19$ &   $2$ &  $114$ & $\neq 1$ & $3$ & $\beta$  & $(3)$   & $(3,3)$    \\
   $95$ &  $19$ &   $5$ &  $285$ & $\neq 1$ & $3$ & $\beta$  & $(3)$   & $(3,3)$    \\
  $146$ &  $73$ &   $2$ &  $438$ & $\neq 1$ & $3$ & $\beta$  & $(3)$   & $(3,3)$    \\
  $209$ &  $19$ &  $11$ &  $627$ & $=1$     & $1$ & $\alpha$ & $(9,3)$ & $(9,9,3)$  \\
 $1577$ &  $19$ &  $83$ & $4731$ & $=1$     & $1$ & $\alpha$ & $(9,3)$ & $(9,9,3)$  \\
 $2147$ &  $19$ & $113$ & $6441$ & $=1$     & $3$ & $\beta$ & $(9,3)$ & $(27,9,3)$  \\
 $3287$ &  $19$ & $173$ & $9861$ & $\neq 1$ & $3$ & $\beta$  & $(3)$   & $(3,3)$    \\
\hline
\end{tabular}
\end{center}
\end{table}

\begin{remark}
\label{rem3}
In Table \ref{tb3}, if the $3$-class group $Cl_{3}(N)$ is of type $(3,3)$, then $N^{\ast}= N_3^{(1)}$.
In this case, the cubic residue symbol $\left(\dfrac{q_{1}}{p_1}\right)_3$ is different from $1$
if and only if $3$ divides exactly the class number of $L$, which implies that $u=3$.
Consequently, Table \ref{tb3} confirms the previous results by Ismaili and El Mesaoudi
\cite[Theorem 1, case $(iv)$, p. 157; Theorem 5 and Corollary 4, p. 165]{IsEM}.
\end{remark}

%%%%%%%%%%%%%%%%%%%%%%%%%%%%%%%%%%%%%%%%%%%%%%%%%%%%%

%%%%%%%%%%%%%%%%%%%%%%%%%%%%%%%%%%%%%%%%%%%%%%%%%%%%%

%\newpage
%\paragraph{}
%%%%%%%%%%%%%%%%%%%%%%%%%%%%%%%%%%%%%%%%%%%%%%%%%%%%%

\noindent
Table \ref{tb4} treats conductors $f=9p_1q_1$ of species 1a, where $p_1\equiv 4,7\,(\mathrm{mod}\,9)$ and $q_1\equiv -1\,(\mathrm{mod}\,3)$.
See Theorem
\ref{thm:Ismaili1},
item (5), and the fourth line of Equation
\eqref{eqn:Rank2}.

\renewcommand{\arraystretch}{1.0}

\begin{table}[ht]
\caption{Radicands $D=3^ep_1q_1$ of species 1a with $m=4$}
\label{tb4}
\begin{center}
\begin{tabular}{|rlrrr|ccccc|cc|}
\hline
 $D$    & $3^e$ & $p_1$ & $q_1$ & $f$    & $\left(\dfrac{3}{p_1}\right)_3$ & $\left(\dfrac{q_{1}}{p_1}\right)_3$ & $\left(\dfrac{3^{e}q_{1}}{p_1}\right)_3$ & $u$ & PFT & $Cl_{3}(L)$ & $Cl_{3}(N)$   \\
\hline
  $342$ & $3^2$ &  $19$ &   $2$ &   $342$ & $\neq 1$ & $\neq 1$ & $=1$     & $1$ & $\alpha$ & $(9,3)$  & $(9,9,3)$   \\
  $555$ & $3$   &  $37$ &   $5$ &  $1665$ & $\neq 1$ & $\neq 1$ & $\neq 1$ & $3$ & $\beta$  & $(3)$    & $(3,3)$     \\
 $3663$ & $3^2$ &  $37$ &  $11$ &  $3663$ & $\neq 1$ & $=1$     & $\neq 1$ & $3$ & $\beta$  & $(3)$    & $(3,3)$     \\
 $5757$ & $3$   &  $19$ & $101$ & $17271$ & $\neq 1$ & $\neq 1$ & $=1$     & $3$ & $\beta$  & $(3,3)$  & $(9,3,3)$   \\
 $6441$ & $3$   &  $19$ & $113$ & $19323$ & $\neq 1$ & $=1$     & $\neq 1$ & $3$ & $\beta$  & $(3)$    & $(3,3)$     \\
 $9519$ & $3$   &  $19$ & $167$ & $28557$ & $\neq 1$ & $\neq 1$ & $=1$     & $1$ & $\alpha$ & $(27,3)$ & $(27,27,3)$ \\
\hline
\end{tabular}
\end{center}
\end{table}

\begin{remark}
\label{rem4}
In Table \ref{tb4}, if the $3$-class group $Cl_{3}(N)$ is of type $(3,3)$, then $N^{\ast}= N_3^{(1)}$.
In this case, the cubic residue symbol $\left(\dfrac{3^{e}q_{1}}{p_1}\right)_3$ is different from $1$
if and only if $3$ divides exactly the class number of $L$, which implies that $u=3$.
Consequently, Table \ref{tb4} confirms the previous results by Ismaili and El Mesaoudi
\cite[Theorem 1, case $(iv)$, p. 157; Theorem 8 and Corollary 7, p. 174]{IsEM}.
\end{remark}

%\paragraph{}
%%%%%%%%%%%%%%%%%%%%%%%%%%%%%%%%%%%%%%%%%%%%%%%%%%%%%

Table \ref{tb5} deals with conductors $f=p_1q_1q_2$ of species 2, where $p_1\equiv 1\,(\mathrm{mod}\,9)$ and $q_1,q_2\equiv 2,5\,(\mathrm{mod}\,9)$.
See Theorem
\ref{thm:Ismaili1},
item (6), and the fifth line of Equation
\eqref{eqn:Rank2}.
Other examples are given in
\cite[Exm. 7.1, p. 113]{AMI},
namely the M0-fields \(D=12673=19\cdot 23\cdot 29\), \(20539=19\cdot 23\cdot 47\).

\renewcommand{\arraystretch}{1.0}
 
\begin{table}[ht]
\caption{Radicands $D=p_1q_1q_2$ of species 2 with $m=2$}
\label{tb5}
\begin{center}
\begin{tabular}{|rrrrr|ccccc|cc|}
\hline
 $D$     & $p_1$ & $q_1$ & $q_2$ & $f$     & $\left(\dfrac{q_{1}}{p_1}\right)_3$ & $\left(\dfrac{q_{2}}{p_1}\right)_3$ & $\left(\dfrac{q_{1}q_{2}}{p_1}\right)_3$ & $u$ & PFT & $Cl_{3}(L)$ & $Cl_{3}(N)$   \\
\hline
   $190$ &  $19$ &   $2$ &   $5$ &   $190$ & $\neq 1$ & $\neq 1$ & $\neq 1$ & $3$ & $\beta$  & $(3)$    & $(3,3)$      \\
   $874$ &  $19$ &   $2$ &  $23$ &   $874$ & $\neq 1$ & $\neq 1$ & $=1$     & $3$ & $\beta$  & $(3,3)$  & $(9,3,3)$    \\
  $1558$ &  $19$ &   $2$ &  $41$ &  $1558$ & $\neq 1$ & $\neq 1$ & $\neq 1$ & $3$ & $\beta$  & $(3)$    & $(3,3)$      \\
  $9361$ &  $37$ &  $11$ &  $23$ &  $9361$ & $=1$     & $=1$     & $=1$     & $1$ & $\alpha$ & $(27,3)$ & $(27,9,3,3)$ \\
 $24679$ &  $37$ &  $23$ &  $29$ & $24679$ & $=1$     & $=1$     & $=1$     & $3$ & $\beta$  & $(9,3)$  & $(9,9,3,3)$  \\
 $43993$ &  $37$ &  $29$ &  $41$ & $43993$ & $=1$     & $\neq 1$ & $\neq 1$ & $3$ & $\beta$  & $(3)$    & $(3,3)$      \\
\hline
\end{tabular}
\end{center}
\end{table}

\begin{remark}
\label{rem5}
In Table \ref{tb5}, if the $3$-class group $Cl_{3}(N)$ is of type $(3,3)$, then $N^{\ast}= N_3^{(1)}$.
In this case, the cubic residue symbol $\left(\dfrac{q_{1}q_{2}}{p_1}\right)_3$ is different from $1$
if and only if $3$ divide exactly the class number of $L$, and imply that $u=3$.
Consequently, Table \ref{tb5} confirms the previous results by Ismaili and El Mesaoudi
\cite[Theorem 1, case $(iii)$, p. 157; Theorem 4 and Corollary 3, pp. 164--165]{IsEM}.
\end{remark}

%\noindent
An example for Theorem
\ref{thm:Ismaili1},
item (8), and the seventh line of Equation
\eqref{eqn:Rank2}
is given by the M0-field \(D=52417=23\cdot 43\cdot 53\) in
\cite[Exm. 7.3, p. 115]{AMI},
which belongs to a doublet \((m=2)\).

%%%%%%%%%%%%%%%%%%%%%%%%%%%%%%%%%%%%%%%%%%%%%%%%%%%%%

\paragraph{}
In the outlook of the paper
\cite[\S 6, Outlook, Figure 3, Scenario III, p. 273]{AMITA},
we presented a scenario for the relative $3$-genus field $N^{\ast}$
in the case where the $3$-class group $Cl_{3}(N)$ is of type $(3,3)$,
and which is relevant for remarks
\ref{rem1}, \ref{rem2}, \ref{rem3}, \ref{rem4} and \ref{rem5}.

\begin{remark}
\label{rmk:Ismaili3}
It should be pointed out that Ismaili \cite{Is} has investigated
another related scenario for the relative $3$-genus field $N^{\ast}$.
If the conductor $f$ is divisible by exactly one prime $p\equiv 1\,(\mathrm{mod}\,3)$, that is $s=1$,
but is not contained in Theorem \ref{thm:Ismaili1},
and if $C_{k,3}\simeq (3,3)$ is elementary bicyclic,
then the genus field coincides with the Hilbert $3$-class field, $N^{\ast}=N_3^{(1)}$,
and the composita $N\cdot L_3^{(1)}=N\cdot (L^{\sigma})_3^{(1)}=N\cdot (L^{\sigma^2})_3^{(1)}=K_4$
coincide with one of the four unramified cyclic cubic extensions $K_1,\ldots,K_4$ of $N$ within $N_3^{(1)}$,
as illustrated in Figure \ref{fig:IsmailiType3}
and studied in detail by Ismaili and El Mesaoudi \cite{IsEM}.

\end{remark}
 
%\newpage

%--------------------------------------------------------------------------------

\begin{figure}[ht]
\caption{Case where  $N^{\ast}=N_3^{(1)}$}
\label{fig:IsmailiType3}

{\tiny

\setlength{\unitlength}{1.0cm}
\begin{picture}(15,11)(-10.5,-10)

% scale of degrees
\put(-10,0.5){\makebox(0,0)[cb]{Degree}}

\put(-10,-2){\vector(0,1){2}}

\put(-10,-2){\line(0,-1){7}}
\multiput(-10.1,-2)(0,-1){8}{\line(1,0){0.2}}

%***********************************************
\put(-10.2,-2){\makebox(0,0)[rc]{\(54\)}}
\put(-9.8,-2){\makebox(0,0)[lc]{Hilbert $3$-class field and relative genus field of the normal closure}}
\put(-10.2,-4){\makebox(0,0)[rc]{\(18\)}}
\put(-9.8,-3.8){\makebox(0,0)[lc]{unramified}}
\put(-9.8,-4.0){\makebox(0,0)[lc]{cyclic cubic}}
\put(-9.8,-4.2){\makebox(0,0)[lc]{extensions of $N$}}
\put(-10.2,-5){\makebox(0,0)[rc]{\(9\)}}
\put(-9.8,-4.9){\makebox(0,0)[lc]{Hilbert $3$-class fields}}
\put(-9.8,-5.2){\makebox(0,0)[lc]{of \(L,L^{\sigma},L^{\sigma^2}\)}}
\put(-10.2,-6){\makebox(0,0)[rc]{\(6\)}}
\put(-9.8,-6){\makebox(0,0)[lc]{normal closure of $L$}}
\put(-10.2,-7){\makebox(0,0)[rc]{\(3\)}}
\put(-9.8,-7){\makebox(0,0)[lc]{conjugate pure cubic fields}}
\put(-10.2,-8){\makebox(0,0)[rc]{\(2\)}}
\put(-9.8,-8){\makebox(0,0)[lc]{cyclotomic quadratic field}}
\put(-10.2,-9){\makebox(0,0)[rc]{\(1\)}}
\put(-9.8,-9){\makebox(0,0)[lc]{rational base field}}
%***********************************************

%{\normalsize
%\put(-3,0){\makebox(0,0)[cc]{Scenario III: \quad \(k^{\ast}=k_1\)}}
%}

%***********************************************

% the cyclotomic tower

% degree 54 vertex
\put(0,-2){\circle*{0.2}}
\put(0.2,-2){\makebox(0,0)[lc]{\(N_3^{(1)}=N^{\ast}\)}}

% directed edges from left to right
\put(0,-2){\line(-3,-2){3}}
\put(0,-2){\line(-1,-2){1}}
\put(0,-2){\line(1,-2){1}}
\put(0,-2){\line(3,-2){3}}

% degree 18 vertices from left to right
\multiput(-3,-4)(2,0){4}{\circle*{0.2}}
\put(-3.2,-3.9){\makebox(0,0)[rc]{\(L_3^{(1)}\cdot N=(L^{\sigma})_3^{(1)}\cdot N=(L^{\sigma^2})_3^{(1)}\cdot N\)}}
\put(-2.8,-4){\makebox(0,0)[lc]{\(=K_4\)}}
\put(-1.2,-4){\makebox(0,0)[rc]{\(K_3\)}}
\put(1.2,-4){\makebox(0,0)[lc]{\(K_2\)}}
\put(3.2,-4){\makebox(0,0)[lc]{\(K_1\)}}

% directed edges from left to right
\put(0,-6){\line(-3,2){3}}
\put(0,-6){\line(-1,2){1}}
\put(0,-6){\line(1,2){1}}
\put(0,-6){\line(3,2){3}}

% degree 6 vertex
\put(0,-6){\circle*{0.2}}
\put(0.2,-6){\makebox(0,0)[lc]{\(N\)}}

% directed edge
\put(0,-6){\line(0,-1){2}}

% degree 2 vertex
\put(0,-8){\circle*{0.2}}
\put(0.2,-8){\makebox(0,0)[lc]{\(K\)}}

%***********************************************
% connecting directed edges 
\put(-5,-5){\line(2,1){2}}
\put(-2,-7){\line(2,1){2}}
\put(-2,-9){\line(2,1){2}}
%***********************************************

% the rational tower

% conjugate degree 9 vertices
\put(-5,-5){\circle{0.2}}
\put(-5.2,-5.1){\makebox(0,0)[rc]{\(L_3^{(1)}\)}}
\put(-4.8,-5){\makebox(0,0)[lc]{\((L^{\sigma})_3^{(1)},(L^{\sigma^2})_3^{(1)}\)}}

% directed edge
\put(-2,-7){\line(-3,2){3}}

% conjugate degree 3 vertices
\put(-2,-7){\circle{0.2}}
\put(-2.2,-7.1){\makebox(0,0)[rc]{\(L\)}}
\put(-1.8,-7){\makebox(0,0)[lc]{\(L^{\sigma},L^{\sigma^2}\)}}

% directed edge
\put(-2,-7){\line(0,-1){2}}

% degree 1 vertex
\put(-2,-9){\circle*{0.2}}
\put(-2.2,-9){\makebox(0,0)[rc]{\(\mathbb{Q}\)}}

%***********************************************

\end{picture}

}

\end{figure}

%%%%%%%%%%%%%%%%%%%%%%%%%%%%%%%%%%%%%%%%%%%%

\subsection{Conductors without splitting prime divisors}
\label{ss:NonSplit}

Now we give more details concerning the trailing five lines
of equation \eqref{eqn:Rank2} in our Theorem \ref{thm:Rank2},
where $D$ is only divisible by  primes $q_{j}\equiv -1\,(\mathrm{mod}\,3)$
which do not split in $K$.

\begin{theorem}
\label{thm:Ismaili2}
Let the conductor of $N/K$ be $f=3^\varepsilon\cdot p_1\cdots p_w\cdot q_1\cdots q_J$
as in Equation
\eqref{eqn:PrimeDecCond}
with $0\le\varepsilon\le 2$, $T=w+J\ge 1$,
and pairwise distinct primes
$p_{i}\equiv 1\,(\mathrm{mod}\,3)$ for $1\le i\le w$,
and $q_{i}\equiv -1\,(\mathrm{mod}\,3)$ for $1\le i\le J$.
Briefly denote the multiplicity of $f$ by $m:=m(f)$.
Assume that $w=0$.
Then, $\operatorname{rank}\,(Cl_{3}^{(\sigma)}(N))=2$
$\Longleftrightarrow$ $L$ belongs to one of the following multiplets:

\begin{enumerate}
\item[$(1)$]
quartets, $m=4$, of type $(\beta^x,\gamma^y)$, $x+y=4$,
such that $f=9q_1q_2$ with $q_1\equiv q_2\equiv 8\,(\mathrm{mod}\,9)$,
\item[$(2)$]
quartets, $m=4$, of type $(\beta^x,\gamma^y)$, $x+y=4$,
such that $f=q_1q_2q_3$ with $q_1\equiv q_2\equiv q_3\equiv 8\,(\mathrm{mod}\,9)$,
\item[$(3)$]
triplets, $m=3$, of type $(\beta,\beta,\beta)$   
such that $f=3q_1q_2q_3$ with $q_1,q_2,q_3\equiv 2,5\,(\mathrm{mod}\,9)$,
\item[$(4)$]
doublets, $m=2$, of type $(\beta,\beta)$   
such that $f=3q_1q_2q_3$ with $q_1,q_2\equiv 2,5\,(\mathrm{mod}\,9)$ and $q_3\equiv 8\,(\mathrm{mod}\,9)$,
\item[$(5)$]
quartets, $m=4$, of type $(\beta,\beta,\beta,\beta)$   
such that $f=3q_1q_2q_3$ with $q_1\equiv 2,5\,(\mathrm{mod}\,9)$ and $q_2\equiv q_3\equiv 8\,(\mathrm{mod}\,9)$,
\item[$(6)$]
triplets, $m=3$, of type $(\beta,\beta,\beta)$ such that $f=q_1q_2q_3q_4$ with
$q_1,q_2,q_3,q_4\equiv 2,5\,(\mathrm{mod}\,9)$,
\item[$(7)$]
doublets, $m=2$, of type $(\beta,\beta)$ such that $f=q_1q_2q_3q_4$ with
$q_1,q_2,q_3\equiv 2,5\,(\mathrm{mod}\,9)$ and $q_4\equiv 8\,(\mathrm{mod}\,9)$,
\item[$(8)$]
quartets, $m=4$, of type $(\beta,\beta,\beta,\beta)$ such that $f=q_1q_2q_3q_4$ with
$q_1,q_2\equiv 2,5\,(\mathrm{mod}\,9)$ and $q_3\equiv q_4\equiv 8\,(\mathrm{mod}\,9)$,
\item[$(9)$]
octets, $m=8$, of type $(\beta,\beta,\beta,\beta,\beta,\beta,\beta,\beta)$ 
such that $f=9q_1q_2q_3$ with $q_{i}\not\equiv 8\,(\mathrm{mod}\,9)$ for some $i\in\lbrace 1,2,3\rbrace$.
\end{enumerate}

\noindent
There exist infinitely many multiplets with conductors of all these shapes $(1)$--$(9)$.
\end{theorem}

%--------------------------------------------------------------------------------

\begin{proof}
The multiplicity $m(f)$ of each conductor is calculated by means of Formula
\eqref{eqn:Multiplicity},
using the sequence $(X_k)_{k\ge -1}=(\frac{1}{2},0,1,1,3,5,11,\ldots)$:
\begin{enumerate}
\item[$(1)$]
For $D=3^{e}q_1^{f_1}q_2^{f_2}\not\equiv\pm 1\,(\mathrm{mod}\,9)$
with $q_1\equiv q_2\equiv -1\,(\mathrm{mod}\,9)$,
we have $f=3^2q_1q_2$ of species 1a.
We must take into consideration that $G=2$, $B=0$, $T=2$,
and we obtain $m(f)=2^T=2^2=4$, a quartet, independently of $G$ and $B$. \\
Here $e,f_1,$ and $f_2 \in   \{1,2\}$.  Since
$\mathbb{Q}(\sqrt[3]{ab^2})=\mathbb{Q}(\sqrt[3]{a^2b})$ for $  a,
b \in \mathbb{N}$, we can choose $e=1$, so $D= 3q_1^{f_1}q_2^{f_2}$. The prime $q_{i}$ is inert in $K$ for each $i \in \lbrace 1,2 \rbrace$,
     and $q_{i}$ is ramified in $L$. The fact that $D \not \equiv \pm1 \pmod  9$ implies that $3$
     is ramified in $L$, then $\lambda$ is ramified in $N/K$.
     Hence $t=3$. As $3=-\zeta_3^2\lambda^2$, then $N=K(\sqrt[3]{x})$
     where $x=\zeta_3^2\lambda^2 \pi_1^{f_1}\pi_2^{f_2}$, and for  $i \in \lbrace 1,2 \rbrace$, $-q_{i}=\pi_{i}$
     is a prime number of $K$. If \(q_1 \equiv q_2 \equiv  -1 \pmod  9\),
      then  all primes $\pi_1$, $\pi_2$ are congruent to $ 1 ~(\bmod~\lambda^3)$,
      so by  \cite[Lem. 3.3, p. 264]{AMITA} we have $\zeta_3$
      is  norm of an element of $N-\lbrace 0\rbrace$ and $q^{*}=1$.
      Thus $\operatorname{rank}\, (Cl_{3}^{(\sigma)}(N))=2$. The opposite implication of the rank condition $\operatorname{rank}\,(Cl_{3}^{(\sigma)}(N))=2$
to this shape of conductor will be proved in \S\ \ref{02}.
%%%%%%%%
\item[$(2)$]
For $D=q_1^{f_1}q_2^{f_2}q_3^{f_3}$
with $q_1\equiv q_2\equiv q_3\equiv 8\,(\mathrm{mod}\,9)$,
we have $f=q_1q_2q_3$ of species 2.
We must take into consideration that $G=3$, $B=0$,
and we obtain $m(f)=2^G\cdot X_{B-1}=2^3\cdot\frac{1}{2}=4$, a quartet. \\
 The prime $q_{i}$ is inert in $K$ for each $i \in   \{1,2,3\}$, because $ q_{i} \equiv 2 \pmod 3$, and $q_{i}$ is ramified in $L$. Since $D \equiv \pm1 \pmod  9$, then $3$ is not ramified in $L$, and then $\lambda$ is not ramified in $N/K$.
   So $t=3$. Since $q_{i} \equiv -1  \pmod 9$, then $\pi_{i}  \equiv 1 ~(\bmod~3\mathcal{O}_{K})$,
   where $-q_{i}=\pi_{i}$ is a prime number of $K$. Thus
 $N=K(\sqrt[3]{x})$ with $x=\pi_1^{f_1} \pi_2^{f_2}\pi_3^{f_3}$.  The fact that all primes $\pi_1$, $\pi_2$ and $\pi_3$ are congruent to $ 1 ~(\bmod~\lambda^3)$ implies by  \cite[Lem. 3.3, p. 264]{AMITA} that $\zeta_3$ is a norm of an element of $N-\lbrace 0\rbrace$ and $q^{*}=1$.
  Hence, $\operatorname{rank}\,(Cl_{3}^{(\sigma)}(N))=2$.
 The opposite implication of the rank condition $\operatorname{rank}\,(Cl_{3}^{(\sigma)}(N))=2$
to this shape of conductor will be proved in \S\ \ref{04}.
%%%%%%%%
\item[$(3)$]
For $D=q_1^{f_1}q_2^{f_2}q_3^{f_3}\not\equiv\pm 1\,(\mathrm{mod}\,9)$
with $q_1,q_2,q_3\equiv 2,5\,(\mathrm{mod}\,9)$,
we have $f=3q_1q_2q_3$ of species 1b.
We must take into consideration that $G=0$, $B=3$,
and we obtain $m(f)=2^G\cdot X_{B}=1\cdot 3=3$, a triplet. \\
 For the equivalence of the rank condition $\operatorname{rank}\,(Cl_{3}^{(\sigma)}(N))=2$
to this shape of conductor, we reason as in case (4) and we will  prove the opposite implication in \S\ \ref{04}.
%%%%%%%%
\item[$(4)$]
For $D=q_1^{f_1}q_2^{f_2}q_3^{f_3}\not\equiv\pm 1\,(\mathrm{mod}\,9)$
such that $(q_1,q_2\equiv 2,5\,(\mathrm{mod}\,9)$ and $q_3\equiv 8\,(\mathrm{mod}\,9))$,
we have $f=3q_1q_2q_3$ of species 1b.
We must take into consideration that $G=1$ and $B=2$,
and we obtain $m(f)=2^G\cdot X_{B}=2\cdot 1=2$, a doublet. \\
For each $i \in   \{1,2,3\}$, $q_{i}$ is inert in $K$ because $ q_{i} \equiv 2 \pmod 3$, and $q_{i}$ is ramified in $L$. As $D \not \equiv \pm1 \pmod  9$, $3$ is  ramified in $L$, so $\lambda$ is  ramified in $N/K$. Then $t=4$.\\ 
   If  $\exists i \in   \{1,2,3\}$ such that  \(q_{i} \not \equiv -1  \pmod 9\), then we have $-q_{i}=\pi_{i}$ is a prime number of $K$. Put $x=\pi_1^{f_1} \pi_2^{f_2}\pi_3^{f_3}$, then
 $N=K(\sqrt[3]{x})$.  According to  \cite[Lem. 3.2, p. 264]{AMITA},  there exists a  prime $\pi_{i}$ not congruent to $ 1 ~(\bmod~\lambda^3)$, then by  \cite[Lem. 3.3, p. 264]{AMITA}, $\zeta_3$ is not a norm of an element of $N-\lbrace 0\rbrace$ and $q^{*}=0$. We conclude that $\operatorname{rank}\,(Cl_{3}^{(\sigma)}(N))=2$.
 The opposite implication of the rank condition $\operatorname{rank}\,(Cl_{3}^{(\sigma)}(N))=2$
to this shape of conductor will be proved in \S\ \ref{04}.
%%%%%%%%
\item[$(5)$]
For $D=q_1^{f_1}q_2^{f_2}q_3^{f_3}\not\equiv\pm 1\,(\mathrm{mod}\,9)$
such that $(q_1\equiv 2,5\,(\mathrm{mod}\,9)$ and $q_2,q_3\equiv 8\,(\mathrm{mod}\,9))$,
we have $f=3q_1q_2q_3$ of species 1b.
We must take into consideration that $G=2$ and $B=1$,
and in the two cases we obtain $m(f)=2^G\cdot X_{B}=4\cdot 1=4$, a quartet.
%%%%%%%%
\item[$(6)$]
For $D=q_1^{f_1}q_2^{f_2}q_3^{f_3}q_4^{f_4}\equiv\pm 1\,(\mathrm{mod}\,9)$
such that $q_1,q_2,q_3,q_4\equiv 2,5\,(\mathrm{mod}\,9)$,
we have $f=q_1q_2q_3q_4$ of species 2.
We must take into consideration that $G=0$, $B=4$,
and we obtain $m(f)=2^G\cdot X_{B-1}=1\cdot 3=3$, a triplet. \\
For the equivalence of the rank condition $\operatorname{rank}\,(Cl_{3}^{(\sigma)}(N))=2$
to this shape of conductor, we reason as in case (8) and we will  prove the opposite implication in \S\ \ref{07}.
%%%%%%%%
\item[$(7)$]
For $D=q_1^{f_1}q_2^{f_2}q_3^{f_3}q_4^{f_4}\equiv\pm 1\,(\mathrm{mod}\,9)$
such that $q_1,q_2,q_3\equiv 2,5\,(\mathrm{mod}\,9)$ and $q_4\equiv 8\,(\mathrm{mod}\,9)$,
we have $f=q_1q_2q_3q_4$ of species 2.
We must take into consideration that $G=1$, $B=3$,
and we obtain $m(f)=2^G\cdot X_{B-1}=2\cdot 1=2$, a doublet. \\
 For the equivalence of the rank condition $\operatorname{rank}\,(Cl_{3}^{(\sigma)}(N))=2$
to this shape of conductor, we reason as in case (6) and we will  prove the opposite implication in \S\ \ref{07}.
%%%%%%%%
\item[$(8)$]
For $D=q_1^{f_1}q_2^{f_2}q_3^{f_3}q_4^{f_4}\equiv\pm 1\,(\mathrm{mod}\,9)$
with $q_1,q_2 \equiv 2,5\,(\mathrm{mod}\,9)$ and $q_3\equiv q_4\equiv 8\,(\mathrm{mod}\,9)$,
we have $f=q_1q_2q_3q_4$ of species 2.
We must take into consideration that $G=2$, $B=2$,
and we obtain $m(f)=2^G\cdot X_{B-1}=2^2\cdot 1=4$, a quartet. \\
For each $i\in\{1,2,3,4\}$, $q_{i}$ is inert in $K$ because $ q_{i} \equiv 2 \pmod 3$, and $q_{i}$ is ramified in $L$. As $D \equiv \pm1 \pmod  9$, $3$ is not ramified in $L$, so $\lambda$ is not ramified in $N/K$.
   Then $t=4$. Put $x=\pi_1^{f_1} \pi_2^{f_2}\pi_3^{f_3}\pi_4^{f_4}$, then
 $N=K(\sqrt[3]{x})$.  If there exists $i \in \{1,2,3,4\} $ such that $ q_{i} \not \equiv -1  \pmod 9,$ then by  \cite[Lem. 3.2, p. 264]{AMITA}, $\pi_{i}  \not \equiv 1 ~(\bmod~3\mathcal{O}_{K})$,
   where $-q_{i}=\pi_{i}$ is a prime number of $K$, then according to  \cite[Lem. 3.3, p. 264]{AMITA}, $\zeta_3$ is not a norm of an element of $N-\lbrace 0\rbrace$ and $q^{*}=0$.
  We conclude that $\operatorname{rank}\,(Cl_{3}^{(\sigma)}(N))=2$. The opposite implication of the rank condition $\operatorname{rank}\,(Cl_{3}^{(\sigma)}(N))=2$
to this shape of conductor will be proved in \S\ \ref{07}.
%%%%%%%%
\item[$(9)$]
For $D=3^{e}q_1^{f_1}q_2^{f_2}q_3^{f_3}\not\equiv\pm 1\,(\mathrm{mod}\,9)$
such that $\exists i\in\lbrace 1,2,3\rbrace\mid q_{i}\not\equiv 8\,(\mathrm{mod}\,9)$,
we have $f=9q_1q_2q_3$ of species 1a, with $T=3$, 
and we obtain $m(f)=2^T=2^3=8$, an octet, independently of $G$ and $B$.\\
 We reason as in case (3), we get $t=4$ and $q^{*}=0$, then $\operatorname{rank}\,(Cl_{3}^{(\sigma)}(N))=2$.
 The opposite implication of the rank condition $\operatorname{rank}\,(Cl_{3}^{(\sigma)}(N))=2$
to this shape of conductor will be proved in \S\ \ref{04}.
\end{enumerate}

The principal factorization type
is a consequence of the estimates in \cite[Thm. 2.1, p. 254]{AMITA}.
Since $s=0$, we have $0\le R\le 0$ and type $\alpha$ with relative PF is generally forbidden.
For all cases, we have $3\le t\le 4$ and thus $1\le A\le\min(2,t-s)=2$
with $3=3-0\le t-s\le 4-0=4$,
which enables type $\beta$ with absolute PF.
For items $(3)$--$(8)$, there exists a prime factor $q\equiv 2,5\,(\mathrm{mod}\,9)$,
and type $\gamma$ is impossible,
because $\zeta_3$ can be norm of a unit in $N$ only if the prime factors of $f$
are $3$ or $\ell_j\equiv 1,8\,(\mathrm{mod}\,9)$.
The latter condition is satisfied by items $(1)$--$(2)$, whence type $\gamma$ may occur.

All claims on the infinitude of the various sets of conductors $f$
are consequences of Dirichlet's theorem on prime numbers arising from invertible residue classes.
\end{proof}

%%%%%%%%%%%%%%%%%%%%%%%%%%%%%%%%%%%%%%%%

%\newpage
%-----------------------------------------------------------------------

\noindent
\textbf{Examples 2.} Several M0-fields
\cite[Dfn. 4.3, p. 105]{AMI}
have radicands $D$ of the shape in the eleventh line of Equation
\eqref{eqn:Rank2}:
\(9922=2\cdot 11^2\cdot 41\), \(38686=2\cdot 23\cdot 29^2\)
belong to Theorem
\ref{thm:Ismaili2},
item (3),
\(850=2\cdot 5^2\cdot 17\), \(6358=2\cdot 11\cdot 17^2\),
\(17986=2\cdot 17\cdot 23^2\), \(94162=2\cdot 23^2\cdot 89\)
belong to Theorem
\ref{thm:Ismaili2},
item (4),
and \(61268=2^2\cdot 17^2\cdot 53\) belongs to Theorem
\ref{thm:Ismaili2},
item (5).
The M0-field \(D=55522=2\cdot 17\cdot 23\cdot 71\) realizes Theorem
\ref{thm:Ismaili2},
item (8) and line thirteen of Equation
\eqref{eqn:Rank2}.
So Equation
\eqref{eqn:Rank2}
indeed gives rise to numerous M0-fields
\cite{AMI}.
See sequence A363699 in OEIS
\cite{OEIS2023}.

%%%%%%%%%%%%%%%%%%%%%%%%%%%%%%%%%%

%-----------------------------------------------------------------------

\section{Proof of the Main Theorem \ref{thm:Rank2}}
\label{s:Proof}
Let $N=\mathbb{Q}(\sqrt[3]{D},\zeta_3)$ be the normal closure of the pure cubic field $L=\mathbb{Q}(\sqrt[3]{D})$,
where $D>1$ is a cube free positive integer, $K=\mathbb{Q}(\zeta_3)$, and
$Cl_{3}(N)$ be the \(3\)-class group of \(N\).
Let $N_3^{(1)}$ be the maximal abelian unramified $3$-extension of $N$.
It is known that $N_3^{(1)}/K$ is Galois, and according to class field theory:

\begin{equation}
\label{IsomCk}
\operatorname{Gal}\left(N_3^{(1)}/N \right)\cong Cl_{3}(N).
\end{equation}  
   
We denote by $N^{\ast}$ the maximal abelian extension of $K$ contained in $N_3^{(1)}$,
which is called the \textit{relative genus field} of $N/K$ (see \cite{Ge1975}, \cite{Ge1976} or \cite{Hz}). 
 
  It is known that the commutator subgroup of $\operatorname{Gal}\left(N_3^{(1)}/K \right)$
  coincides with $\operatorname{Gal}\left(N_3^{(1)}/ N^{\ast} \right)$ and then: 
  $$\operatorname{Gal}\left(N^{\ast}/K\right) \cong \operatorname{Gal}\left(N_3^{(1)}/K\right) \Bigm/ \operatorname{Gal}\left(N_3^{(1)}/ N^{*}\right).  $$
 
Let $\sigma$ be a generator of $\operatorname{Gal}\left(N/K\right)$,
and let $Cl_{3}^{1-\sigma}(N)$ be the subgroup of $Cl_{3}(N)$ defined by
$Cl_{3}^{1-\sigma}(N) = \lbrace \mathcal{A}^{1-\sigma} \mid \mathcal{A}\in Cl_{3}(N) \rbrace$,
which is called the \textit{principal genus} of $Cl_{3}(N)$.
The fact that $N/K$ is abelian and that $N \subseteq N^{\ast}$ implies that  
$\operatorname{Gal}\left(N_3^{(1)}/ N^{\ast} \right)$ coincides with $ Cl_{3}^{1-\sigma}(N)$,
by the aid of the isomorphism \eqref{IsomCk} above, and by Artin's reciprocity law.
See Figure
\ref{fig:GenusField}.
 
%\newpage

%\begin{center}
%Figure: The Galois correspondence between \\ the principal genus $Cl_{3}^{1-\sigma}(N)$ and the relative genus field $N^\ast$: 
%\end{center}
\begin{figure}[ht]
\caption{Galois correspondence between principal genus $Cl_{3}^{1-\sigma}(N)$ and relative genus field $N^\ast$}
\label{fig:GenusField}
\[
  \xymatrix@R=0,5cm{
       {}&{}&{}&{}&{} \\
       {}& N_3^{(1)} \ar@{<-}[d] &{}& \{1\} \ar@{->}[d] &{} \\
       {}& N^{\ast} \ar@{<-}[d]  &{}& \mathrm{Gal(N_3^{(1)}/N^{\ast})}\simeq Cl_{3}^{1-\sigma}(N) \ar@{->}[d] &{} \\
       {}& N                     &{}& \mathrm{Gal(N_3^{(1)}/N)}\simeq Cl_{3}(N) &{} \\
       {}&{}&{}&{}&{} \\
  }
\]
\end{figure}
Then
$\operatorname{Gal}\left( N^{\ast}/ N  \right) \cong Cl_{3}(N)/Cl_{3}^{1-\sigma}(N).$
Let
$Cl_{3}^{(\sigma)}(N) = \lbrace \mathcal{A} \in Cl_{3}(N) \mid \mathcal{A}^{\sigma} = \mathcal{A} \rbrace$
be the $3$-group of \textit{ambiguous ideal classes} of $N/K$.
Since the Sylow $3$-subgroup of the ideal class group of $K$ is reduced to $\lbrace 1\rbrace$, and by the exact sequence:
$$1 \longrightarrow Cl_{3}^{(\sigma)}(N)  \longrightarrow 
Cl_{3}(N) \overset{1-\sigma}{\longrightarrow } Cl_{3}(N) \longrightarrow Cl_{3}(N)/Cl_{3}^{1-\sigma}(N) \longrightarrow 1 $$
we see that $\operatorname{Gal}\left( N^{\ast}/ N\right) \cong Cl_{3}^{(\sigma)}(N)$.

We assume $\operatorname{Gal}\left(N^{\ast}/N\right)\cong\mathbb{Z}/3\mathbb{Z}\times\mathbb{Z}/3\mathbb{Z}$.
In Equation $(3.2)$ of \cite[p. 55]{Ge1975}, the integer $D$ is written in the following form:
\begin{eqnarray}
\label{eq2}
D &=& 3^e\cdot p_1^{e_1}\cdots p_v^{e_v}\cdot p_{v+1}^{e_{v+1}}\cdots p_w^{e_w}\cdot q_1^{f_1}\cdots q_I^{f_I}\cdot q_{I+1}^{f_{I+1}}\cdots q_J^{f_J},
\end{eqnarray}
where $p_i$ and $q_i$ are positive rational primes such that:
\[
\left\{
\begin{array}{ll}
 p_i \equiv 1\,(\mathrm{mod}\,9), & \quad \text{ for } \quad 1\leq i\leq v, \\
 p_i \equiv 4 \text{ or } 7\,(\mathrm{mod}\,9), & \quad \text{ for } \quad v+1\leq i\leq w,\\
 q_i \equiv 8\,(\mathrm{mod}\,9), & \quad \text{ for } \quad 1\leq i\leq I,\\
 q_i \equiv 2 \text{ or } 5\,(\mathrm{mod}\,9), & \quad \text{ for } \quad I+1\leq i\leq J,\\
 e_i = 1 \text{ or } 2, & \quad \text{ for } \quad 1\leq i\leq w,\\
 f_i = 1 \text{ or } 2, & \quad \text{ for } \quad 1\leq i\leq J,\\
 e   = 0, 1 \text{ or } 2.
\end{array}
\right.
\]

\noindent
Since \(\operatorname{rank}\,(Cl_{3}^{(\sigma)}(N))=2\), then Lemma $3.1$ of \cite[p. 55]{Ge1975} gives the following cases:
   \begin{itemize}
   \item \textbf{Case 1:} \(2w+J=2\);
   \item \textbf{Case 2:} \(2w+J=3\);
   \item \textbf{Case 3:} \(2w+J=4\);
   \end{itemize}
where $w$ and $J$ are defined in Equation \eqref{eq2}. We shall treat the three cases above as follows:
    
    \subsection{Case 1:}
%%%%%%%%%%%%%%%%%%%%%%%%%%%%%%%%%%%%%
%%%%%%%%%%%%%%%% CASE 1 %%%%%%%%%%%%%%%%%%%%%%
%%%%%%%%%%%%%%%%%%%%%%%%%%%%%%%%%%%%%%
   In the case $1$, $2w+J=2$, we either have ($w=1$ and $J=0$) treated in \S\ \ref{01}  or ($w=0$ and $J=2$) treated in \S\ \ref{02} .
    \subsubsection{Radicands divisible by one splitting rational prime}
    \label{01}

    If $w=1$ and $J=0$, then
     $$D=3^{e}p^{e_1},$$
      where $p$ is a prime number such that
   \(p \equiv 1 (\bmod 3)\), $e\in   \{0,1,2\}$ and $e_1 \in   \{1,2\}$. Then: 
  
 \paragraph*{Species 2:}
   If \(D \equiv \pm1 \pmod  9\), then we necessary have \(p \equiv 1 \pmod  9\) and $e=0$. So the integer $D$
   can be written in the form  $D=p^{e_1}$, with $ p \equiv 1  \pmod 9$ and  $e_1\in   \{1,2\}$. According to \cite[Thm.1.1, p. 251]{AMITA},  $\operatorname{rank}\,(Cl_{3}^{(\sigma)}(N))=1$, which contradicts our hypothesis.
    
   \paragraph*{Species 1:}
     If \(D \not \equiv \pm1 \pmod  9\), then according to \cite[p. 266]{AMITA}, the integer $D$ is written in one of the following forms:
    \[ d = \left\{
  \begin{array}{l l l}
    p^{e_1} & \quad \text{with \ $ p \equiv 4 ~or~ 7  \pmod 9.$     }\\
   3^{e}p^{e_1} \not \equiv \pm1 \pmod  9 & \quad \text{with \ $ p \equiv 1   \pmod  3. $     }\\

  \end{array} \right.\]
where $e,e_1 \in   \{1,2\}$.  \\

Assume that $D=p^{e_1} $ or $D=3^{e}p^{e_1}$, with $p \equiv 4 ~or~ 7  \pmod 9$, then by \cite[Thm.1.1, p. 251]{AMITA}, we have  $\operatorname{rank}\,(Cl_{3}^{(\sigma)}(N))=1$, which is a contradiction.

 Hence, the possible form of  $D$ in this situation is:
  $$D = 3^{e}p^{e_1} \not \equiv \pm1 \pmod  9 \ \text{with} \    p \equiv  1 \pmod 9.  
    $$
where $e,e_1 \in   \{1,2\}$, which is the first form of $D$ in Theorem \ref{thm:Rank2}.\\

%%%%%%%%%%%%%%%%%%%%%%%%%%%%%%
%%%%%%%%%$w=0$ and $J=2$%%%%%%%%%%%%%%%%%%%%

\subsubsection{Radicands not divisible by a splitting rational prime}
\label{02}
 If $w=0$ and $J=2$, then
  $$D=3^{e}q_1^{f_1}q_2^{f_2},$$
   with
    \(q_1 \equiv q_2 \equiv  2 (\bmod \ 3)\), $e\in   \{0,1,2\}$ and $f_1,f_2 \in   \{1,2\}$, then:
  \paragraph*{Species 2:}
    If \(D \equiv \pm1 \pmod  9\), then according to \cite[p. 266-267]{AMITA}, $e=0$ and \(q_1 \equiv q_2 \equiv -1 \pmod  9\), so  $\operatorname{rank}\,(Cl_{3}^{(\sigma)}(N))=1$ by \cite[Thm. 1.1, p. 251]{AMITA}, which is a contradiction.
     \paragraph*{Species 1:}
    
     If \(D \not \equiv \pm1 \pmod  9\). Assume that 
     $\exists i \in \lbrace 1,2 \rbrace \mid 
     q_{i} \not \equiv   -1 \pmod  9$ such that
      $D=3^{e}q_1^{f_1}q_2^{f_2}\not \equiv \pm1 \pmod  9$, 
     with $e \in   \{0,1,2\}$, $f_1,f_2 \in   \{1,2\}$, then according to \cite[Thm. 1.1, p. 251]{AMITA} we get $\operatorname{rank}\,(Cl_{3}^{(\sigma)}(N))=1$
     which is a contradiction.\\
      Thus, it remains only the form:
      $$D=3^{e}q_1^{f_1}q_2^{f_2} \not\equiv\pm 1\,(\mathrm{mod}\,9),
     \ \ \text{such that} \ q_1 \equiv q_2 \equiv  -1 \pmod  9,$$ 
     where $e,f_1,$ and $f_2 \in   \{1,2\}$, which is the $9^{th}$ form of $D$ in Theorem \ref{thm:Rank2}.

%%%%%%%%%%%%%%%%%%%%%%%%%%%%%%%%%%%%%%%
%%%%%%%%%%%%%%%% CASE 2 %%%%%%%%%%%%%%%%%%%%%%%
%%%%%%%%%%%%%%%%%%%%%%%%%%%%%%%%%%%%%%%

 \subsection{Case 2:}
 In the case $2$, $2w+J=3$, we either have ($w=1$ and $J=1$) treated in \S\ \ref{03}  or ($w=0$ and $J=3$) treated in \S\ \ref{04} .

%%%%%%%%%%%%%%%%%%%%%%%%%%%%%%%%%

\subsubsection{Radicands divisible by one splitting rational prime}
\label{03}
  If \(w=1\) and \(J=1\), then
   $$D=3^{e}p^{e_1}q^{f_1},$$
    where $p$ and $q$ are prime numbers such that  \(p \equiv 1 (\bmod 3) \)
    and \(q\equiv  2 \pmod 3 \), $e\in   \{0,1,2\}$ and $e_1,f_1 \in   \{1,2\}$. Then:
   
\paragraph*{Species 2:}
If \(D \equiv \pm1 \pmod  9 \):
   \begin{itemize}
%%%%---------------------------
   \item[1)] If \(p \equiv 4\, \ \text{or} \ \,7 \pmod  9\) and \(q \equiv -1 (\bmod
   9)\), then according to \cite[p. 267]{AMITA} we get a contradiction.
   %%%%----------------------
  \item[2)] If \(p \equiv -q \equiv 4\, \ \text{or}\, \ 7 \pmod  9\), then according to \cite[p. 267]{AMITA} we have $D= p^{e_1}q^{f_1}\equiv \pm1 \pmod  9$, where $e_1,f_1 \in   \{1,2\}$. But in this case, $\operatorname{rank}\,(Cl_{3}^{(\sigma)}(N))=1$ by \cite[Thm. 1.1, p. 251]{AMITA}. Rejected case.
%%%-------------------------
 \item[3)] If $p \equiv -q \equiv 1  \pmod 9$, we get $D \equiv \pm 3^{e} \pmod 9$, so we have necessary $e=0$.
 Then $D= p^{e_1}q^{f_1}$, where $e_1,f_1 \in   \{1,2\}$,
 which is the second form of $D$ in Theorem \ref{thm:Rank2}. 
%%------------------------
   \item[4)] If \(p \equiv 1 \pmod  9\) and \(q \equiv 2\, \ \text{or}\, 5 \ (\bmod
   9)\), then by \cite[p. 268]{AMITA} we get a contradiction.
\end{itemize}
%%%%%%%%%%%%%%%%%%%%%%%%%%%%%%%%%%%%%%%%%%%%%%
 \paragraph*{Species 1:} If \(D \not \equiv \pm1 \pmod  9\): 
 \\
  According to \cite[\S\ 5, p. 92]{Ge1976}, \( \operatorname{rank}\,(Cl_{3}^{(\sigma)}(N))=t-2+q^*\),
   where $t$ and $q^*$ are defined in the notations.  Since $p\equiv 1 \ (\bmod \ 3$), then $p=\pi_1\pi_2$, where $\pi_1$ and
$\pi_2$ are two primes of $K$ such that
$\pi_2=\pi_1^{\tau}$ and $\pi_1 \equiv \pi_2 \equiv 1 ~(\bmod~3\mathcal{O}_{K})$,
and $p$ is ramified in $L$, then $\pi_1$ and $ \pi_2$ are ramified in $N$.
Also
 $q$ ramifies in $L$ and $-q=\pi$ is a prime in $K$. Since \(D \not \equiv \pm1 \pmod  9\),
then $3$ is ramified in $L$,
and we have $3\mathcal{O}_{K}=(\lambda)^2$. where $\lambda=1-\zeta_3.$
We get $t=4$. \\
%%----------------------------
   If $p \equiv -q \equiv 1  \pmod 9$, 
 then $\pi \equiv \pi_1 \equiv \pi_2 \equiv 1 ~(\bmod~ \lambda^3)$, and according to \cite[Lem. 3.3, p. 264]{AMITA}, $\zeta_3$ is
           norm of an element of $N-\lbrace 0\rbrace$, so $q^{*}=1$.
           Thus, $\operatorname{rank}\,(Cl_{3}^{(\sigma)}(N))=3$, which is a contradiction.\\ 
Hence, the forms of $D$ in this situation are 
 \[ D = \left\{
  \begin{array}{l l l}
   p^{e_1}q^{f_1} \not \equiv \pm1 \pmod  9 & \quad \text{with \ $ p \not \equiv 1  \pmod  9 \ or \ q \not \equiv -1 (\bmod
   9) $ }    \\
  3^{e}p^{e_1}q^{f_1} \not \equiv \pm1 \pmod  9 & \quad \text{with \ $ p \not \equiv 1  \pmod  9 \ or \ q \not \equiv -1 (\bmod
   9) $ }    \\
  \end{array} \right.\] 
  where $e,e_1,f_1 \in   \{1,2\},$
  which are the $3^{th}$ and $5^{th}$ forms of $D$ in Theorem \ref{thm:Rank2}.
 %%%%%%%%%%%%%%%%%%%
 %%%%%%$w=0$ and $J=3$%%%%%%%%%%
 %%%%%%%%%%%%%%%%%%%%
 \subsubsection{Radicands not divisible by a splitting rational prime}
 \label{04}
 If $w=0$ and $J=3$, then
  $$D=3^{e}q_1^{f_1}q_2^{f_2}q_3^{f_3},$$ 
  where $q_{i}$ is a prime number such that
    \(q_{i}\equiv  2  \pmod  3 \), $e\in   \{0,1,2\}$ and $f_{i} \in   \{1,2\}$
    for each $i\in   \{1,2,3\}$. Then:

   %%----------------------------
   
   \paragraph*{Species 2:} If \(D \equiv \pm1 \pmod  9 \):
    \begin{itemize}

   \item[1)] If \(q_1 \equiv 2\, \ \text{or} \ \,5 \pmod  9\) and \(q_2 \equiv q_3 \equiv -1 \ (\bmod \
   9)\), then according to \cite[p. 268]{AMITA}, we get a contradiction. 
  \item[2)] If \(q_1 \equiv q_2 \equiv 2\, \ \text{or} \ \,5 \pmod  9\) and \( q_3 \equiv -1 \pmod
   9\), then according to \cite[p. 268]{AMITA}, we have \(e=0\), so  $D= q_1^{f_1}q_2^{f_2}q_3^{f_3} \equiv \pm1 \pmod  9$, where $f_1,f_2, f_3\in   \{1,2\}$. But in this case we have   $\operatorname{rank}\,(Cl_{3}^{(\sigma)}(N))=1$ by \cite[Thm. 1.1, p. 251]{AMITA} which is  absurd.

 \item[3)] If \(q_1 \equiv q_2 \equiv q_3 \equiv -1  \pmod  9\), then  we have necessary \(e=0\). Then $D=q_1^{f_1}q_2^{f_2}q_3^{f_3}$,
 where $f_1,f_2, f_3\in   \{1,2\}$, which is the $10^{th}$ form of $D$ in Theorem \ref{thm:Rank2}. 
   \item[4)] If \(q_1 \equiv q_2 \equiv q_3 \equiv 2 \, \text{or} \,  5  \pmod
   9\), then  by \cite[p. 268]{AMITA} we have \(e=0\). So, $D= q_1^{f_1}q_2^{f_2}q_3^{f_3}\equiv \pm1 \pmod  9$, where $f_1,f_2, f_3\in   \{1,2\}$, and by \cite[Thm. 1.1, p. 251]{AMITA},
$\operatorname{rank}\,(Cl_{3}^{(\sigma)}(N))=1$, which is  absurd.

\end{itemize}

%%--------------------------
 \paragraph*{Species 1:} If \(D \not \equiv \pm1 \pmod  9\):

For each $i \in   \{1,2,3\}$, $q_{i}$ is inert in $K$ because $ q_{i} \equiv 2 \pmod 3$, and $q_{i}$ is ramified in $L$. As $D \not \equiv \pm1 \pmod  9$, $3$ is  ramified in $L$, so $\lambda$ is  ramified in $N/K$. Then $t=4$.\\ 
    If \(q_1 \equiv q_2 \equiv q_3 \equiv -1  \pmod  9\):  
   \\ For each $i \in   \{1,2,3\}$,  $\pi_{i}  \equiv 1 ~(\bmod~3\mathcal{O}_{K})$, where $-q_{i}=\pi_{i}$ is a prime number of $K$. Put $x=\pi_1^{f_1} \pi_2^{f_2}\pi_3^{f_3}$, then
 $N=K(\sqrt[3]{x})$.  Since all the primes $\pi_1$, $\pi_2$ and $\pi_3$ are congruent to $ 1 ~(\bmod~\lambda^3)$, then according to \cite[Lem. 3.3, p. 264]{AMITA}, $\zeta_3$ is a norm of an element of $N-\lbrace 0\rbrace$ and $q^{*}=1$. We conclude that $\operatorname{rank}\,(Cl_{3}^{(\sigma)}(N))=3$, which is  absurd. \\  
 Thus, it remains only the case where there exist $ i \in   \{1,2,3\}$ such that  \(q_{i} \not \equiv -1  \pmod 9\). 
Hence, the forms of $D$ in this situation are 
 \[ D = \left\{
  \begin{array}{l l l}
   q_1^{f_1}q_2^{f_2}q_3^{f_3} \not \equiv \pm1 \pmod  9, & \quad \text{with \ $ \exists i \in   \{1,2,3\} \mid $  \(q_{i} \not \equiv -1  \pmod 9 $, }    \\
  3^{e}q_1^{f_1}q_2^{f_2}q_3^{f_3} \not \equiv \pm1 \pmod  9, & \quad \text{with \ $ \exists i \in   \{1,2,3\} \mid $  \(q_{i} \not \equiv -1  \pmod 9 $, }    \\
  \end{array} \right.\] 
  where $e,f_1, f_2 $ and $ f_3 \in   \{1,2\}$, which are the $11^{th}$ and $12^{th}$ forms of $D$ in Theorem \ref{thm:Rank2}.

%%%%%%%%%%%%%%%%%%%%%%%%%%%%%%%%%%%%%%%
%%%%%%%%%%%%%%%% CASE 3 %%%%%%%%%%%%%%%%%%%%%%%
%%%%%%%%%%%%%%%%%%%%%%%%%%%%%%%%%%%%%%%
 \subsection{Case 3:}
 In the case $3$, $2w+J=4$ ,we either have ($w=1$ and $J=2$) treated in \S\ \ref{05}, or ($w=2$ and $J=0$) treated in \S\ \ref{06}, or ($w=0$ and $J=4$) treated in \S\ \ref{07}.

 %%%%%%%%%%%%%%%%%
%%W=1 et J=2%%%%%%%%%
%%%%%%%%%%%%%%%%%
\subsubsection{Radicands divisible by one splitting rational prime}
\label{05}
    If \(w=1\) and \(J=2\), then
   $$D=3^{e}p^{e_1}q_1^{f_1}q_2^{f_2},$$ 
   where $p$, $q_1$ and $q_2$ are prime numbers such that  \(p \equiv 1 \pmod 3 \)
    and \(q_1 \equiv q_2 \equiv  2 \pmod 3 \), $e\in   \{0,1,2\}$ and $e_1,f_1, f_2 \in   \{1,2\}$. Then:
    \paragraph*{Species 2:}
    If \(D \equiv \pm1 \pmod  9 \):
   \begin{itemize}
%%%%---------------------------
   \item[1)] If \( p \equiv - q_1 \equiv - q_2 \equiv  1 \pmod 9\),\\
   By \cite[\S\ 5, p. 92]{Ge1976}, \( \operatorname{rank}\,(Cl_{3}^{(\sigma)}(N))=t-2+q^*\),
   where $t$ and $q^*$ are defined in the notations. 
   Since $p\equiv 1 \ (\bmod \ 3$), then  $p=\pi_1\pi_2$, where $\pi_1$ and
$\pi_2$ are two primes of $K$ such that
$\pi_2=\pi_1^{\tau}$ and $\pi_1 \equiv \pi_2 \equiv 1 ~(\bmod~3\mathcal{O}_{K})$,
and the prime $p$ is ramified in $L$, then $\pi_1$ and $ \pi_2$ are ramified in $N$.
For each $i \in   \{1,2\}$, $q_{i}$ is inert in $K$ because $ q_{i} \equiv 2 \pmod 3$, and $q_{i}$ is ramified in $L$. As $D  \equiv \pm1 \pmod  9$, $3$ is not  ramified in $L$, so $\lambda$ is not  ramified in $N/K$. So, $t=4$. As $ q_{1} \equiv q_{2} \equiv - 1 \pmod 9$,  then  for each $i \in   \{1,2\}$, $\pi^{'}_{i}  \equiv 1 ~(\bmod~3\mathcal{O}_{K})$, where $-q_{i}=\pi^{'}_{i}$ is a prime number of $K$. Put $x=-\zeta_3^2\lambda^2\pi_1^{e_1} \pi_2^{e_2}\pi_1^{'f_1}\pi_2^{'f_2}$, then
 $N=K(\sqrt[3]{x})$.  Since all the primes $\pi_1$, $\pi_2$,  $\pi_1^{'}$ and $\pi_2^{'}$ are congruent to $ 1 ~(\bmod~\lambda^3)$, then according to \cite[Lem. 3.3, p. 264]{AMITA}, $\zeta_3$ is a norm of an element of $N-\lbrace 0\rbrace$ and $q^{*}=1$. We conclude that $\operatorname{rank}\,(Cl_{3}^{(\sigma)}(N))=3$, which is  absurd.    
 %%------------------
   \item[2)] If \(p \equiv 1 \pmod 9\) and \(q_1 \equiv q_2 \equiv 2\, \ \text{or}\, \ 5 \pmod 9\), then  $D=3^{e}p^{e_1}q_1^{f_1}q_2^{f_2}\equiv \pm 3^{e},\pm 4\times 3^{e} \ or \ \pm 7 \times 3^{e} \pmod 9 $. Since $D \equiv \pm 1 \pmod 9 $, then $e=0$, so the possible form for the integer $D$ is $D=p^{e_1}q_1^{f_1}q_2^{f_2} \equiv \pm 1 \pmod 9$, 
   which is the $6^{th}$ form of $D$ in Theorem \ref{thm:Rank2}.  
   %%%-------------
   \item[3)] If \(p \equiv 1 \pmod 9\), \(q_1 \equiv -1 \pmod 9\) and \( q_2 \equiv 2\, \ \text{or}\, \ 5 \pmod 9\), then
   $D=3^{e}p^{e_1}q_1^{f_1}q_2^{f_2}\equiv \pm 2\times 3^{e} \ or \ \pm 5 \times 3^{e} \pmod 9 $, so $D \not \equiv \pm 1 \pmod 9$, which is  absurd.
   %%%%----------------------
  \item[4)] If \(p \equiv  4\, \ \text{or}\, \ 7 \pmod  9\) and \(q_1 \equiv q_2 \equiv - 1 \pmod 9\). Then
   $D=3^{e}p^{e_1}q_1^{f_1}q_2^{f_2}\equiv \pm 4\times 3^{e} \ or \ \pm 7 \times 3^{e} \pmod 9 $, so $D \not \equiv \pm 1 \pmod 9$, which is  absurd.
  %%%%-----------------
   \item[5)] If \(p \equiv  4\, \ \text{or}\, \ 7 \pmod  9\) and \(q_1 \equiv q_2 \equiv 2\, \ \text{or}\, \ 5 (\bmod 9)\), then  $D=3^{e}p^{e_1}q_1^{f_1}q_2^{f_2}\equiv \pm 3^{e},\pm 4\times 3^{e} \ or \ \pm 7 \times 3^{e} \pmod 9 $. Since $D \equiv \pm 1 \pmod 9 $, then $e=0$, so the possible form of $D$ is $D=p^{e_1}q_1^{f_1}q_2^{f_2} \equiv \pm 1 \pmod 9$,  which is the $7^{th}$ form of $D$ in Theorem \ref{thm:Rank2}.    
   %%---------------
    \item[6)] If \(p \equiv 4\, \ \text{or}\, \ 7  \pmod 9\), \(q_1 \equiv -1 \pmod 9\) and \( q_2 \equiv 2\, \ \text{or}\, \ 5 \pmod 9\),  then  $D=3^{e}p^{e_1}q_1^{f_1}q_2^{f_2}\equiv  3^{e}, 4\times 3^{e} \ or \  7 \times 3^{e} \pmod 9 $,  since $D \equiv \pm 1 \pmod 9 $, then $e=0$, so the form of $D$ is $D=p^{e_1}q_1^{f_1}q_2^{f_2} \equiv \pm 1 \pmod 9$,  which is the $8^{th}$ form of $D$ in Theorem \ref{thm:Rank2}.    
\end{itemize}
%%%%%%%%%%%%%%%%%%%%%%%%%%%%%%%%%%%%%%%%%%%%%%
 \paragraph*{Species 1:}
 If \(D \not \equiv \pm1 \pmod  9\):  
   The fact that \(\operatorname{rank}\,(Cl_{3}^{(\sigma)}(N))=2\) implies that $t \leq 4.$
   We have \(D=3^{e}p^{e_1}q_1^{f_1}q_2^{f_2}\), with \(p \equiv 1 \pmod
   3 \) and \(q\equiv  2 \pmod 3 \).
   We shall calculate the number of prime ideals which are ramified in $N/K$.
   Since $p\equiv 1 \ (\bmod \ 3$), then  $p=\pi_1\pi_2$, where $\pi_1$ and
$\pi_2$ are two primes of $K$ such that
$\pi_2=\pi_1^{\tau}$ and $\pi_1 \equiv \pi_2 \equiv 1 ~(\bmod~3\mathcal{O}_{K})$,
and the prime $p$ is ramified in $L$, then $\pi_1$ and $ \pi_2$ are ramified in $N$.
Also $q_1$ and $q_2$ are inert in $K$,  $q_1$ and $q_2$  are ramified in $L$. Since \(D \not \equiv \pm1 \pmod  9\),
then $3$ is ramified in $L$,
and we have $3\mathcal{O}_{K}=(\lambda)^2$ where $\lambda=1-\zeta_3.$
Hence, $t=5$ which contradicts the fact that $t \leq 4.$ Rejected case.

%%%%%%%%%%%%%%%%%%%%%
%%%%%%%%%%$w=0$ and $J=4$%%%%%%%%%%%%%
%%%%%%%%%%%%%%%%

%%%%%%%%%%
%%%%%% $w=2$ and $J=2$, %%%%%%%
%%%%%%%%%%%%%%%%

\subsubsection{Radicands divisible by two splitting rational primes}
\label{06}
 If $w=2$ and $J=0$, then
  $$D=3^{e}p_1^{e_1}p_2^{e_2}, $$
   where $p_1$ and $p_2$ are prime numbers such that  \(p_1 \equiv p_2 \equiv 1 \pmod 3 \), $e\in   \{0,1,2\}$ and $e_1, e_2 \in   \{1,2\}$. 
   We shall calculate the number of prime ideals which are ramified in $N/K$.
   Since $p_1 \equiv p_2 \equiv 1 \ (\bmod \ 3$), then  $p_1=\pi_1\pi_2$ and $p_2=\pi_3\pi_4$, where $\pi_1$, $\pi_2$, $\pi_3$ and
$\pi_4$ are primes of $K$ such that
$\pi_2=\pi_1^{\tau}$, $\pi_4=\pi_3^{\tau}$, and $\pi_1 \equiv \pi_2 \equiv \pi_3 \equiv \pi_4 \equiv 1 ~(\bmod~3\mathcal{O}_{K})$,
and the primes $p_1$ and $p_2$ are ramified in $L$, then $\pi_1$, $ \pi_2$, $ \pi_3$, and $ \pi_4$ are ramified in $N$.
\paragraph*{Species 2:}
$D  \equiv \pm1 \pmod  9$, $3$ is  not ramified in $L$, so $\lambda$ is not  ramified in $N/K$. Then $t=4$.\\

Assuming  $ p_1 \equiv p_2 \equiv 1 \pmod 9 $. Since $3=-\zeta_3^2\lambda^2$, then
$N=K(\sqrt[3]{x})$, where $x=-\zeta_3^2\lambda^2\pi_1^{e_1} \pi_2^{e_1}\pi_3^{e_2}\pi_4^{e_2}$.  Since all the primes $\pi_1$, $\pi_2$, $\pi_3$ and $\pi_4$ are congruent to $ 1 ~(\bmod~\lambda^3)$, then according to \cite[Lem. 3.3, p. 264]{AMITA}, $\zeta_3$ is norm of an element of $N-\lbrace 0\rbrace$ and $q^{*}=1$. We conclude that $\operatorname{rank}\,(Cl_{3}^{(\sigma)}(N))=3$, which is an absurd. \\
It remains the case where $ p_1 \not \equiv 1 \pmod 9 $ or $ p_2 \not \equiv 1 \pmod 9 $, this case implies that
 $ p_1^{e_1}p_2^{e_2} \equiv 1, 4$ or $7 \pmod 9$, then $D=3^{e}p_1^{e_1}p_2^{e_2} \equiv 3^{e}, 4\times 3^{e}$ or $7\times 3^{e} \pmod 9 $. Since $D  \equiv \pm1 \pmod  9$, we have necessarily $e=0$, then the possible form of the integer $D$ in this case is: 
$$D=p_1^{e_1}p_2^{e_2} \equiv \pm1 \pmod  9, $$
with $e_1, e_2 \in   \{1,2\}$,  which is the $4^{th}$ form of $D$ in Theorem \ref{thm:Rank2}.  

\paragraph*{Species 1:}
 $D \not \equiv \pm1 \pmod  9$, $3$ is   ramified in $L$, so $\lambda$ is   ramified in $N/K$. Then $t=5$. In this case,  we conclude that $\operatorname{rank}\,(Cl_{3}^{(\sigma)}(N))  \geq 3 $, which is  absurd. 
 
\subsubsection{Radicands not divisible by a splitting rational prime}
\label{07}
 If $w=0$ and $J=4$, then
   $$D=3^{e}q_1^{f_1}q_2^{f_2}q_3^{f_3}q_4^{f_4}, $$
    where $q_{i}$ is a prime number such that
    \(q_{i}\equiv  2  \pmod  3 \), $e\in   \{0,1,2\}$ and $f_{i} \in   \{1,2\}$
    for each $i\in   \{1,2,3,4\}$. Then:

   %%----------------------------
   
   \paragraph*{Species 2:} If \(D \equiv \pm1 \pmod  9 \). Then:
   \begin{itemize}
   \item[1)]  If \(q_1 \equiv q_2 \equiv q_3 \equiv q_4 \equiv -1  \pmod  9\):\\
   For each $i \in   \{1,2,3,4\}$, $q_{i}$ is inert in $K$ because $ q_{i} \equiv 2 \pmod 3$, and $q_{i}$ is ramified in $L$. As $D  \equiv \pm1 \pmod  9$, $3$ is not  ramified in $L$, so $\lambda$ is not ramified in $N/K$. Then $t=4$. Since \(q_1 \equiv q_2 \equiv q_3 \equiv q_4 \equiv -1  \pmod  9\), then for each $i \in   \{1,2,3,4\}$,  $\pi_{i}  \equiv 1 ~(\bmod~3\mathcal{O}_{K})$, where $-q_{i}=\pi_{i}$ is a prime number of $K$. Put $x=-\zeta_3^2\lambda^2\pi_1^{f_1} \pi_2^{f_2}\pi_3^{f_3}\pi_4^{f_4}$, then
 $N=K(\sqrt[3]{x})$.  Since all the primes $\pi_1$, $\pi_2$, $\pi_3$ and $\pi_4$ are congruent to $ 1 ~(\bmod~\lambda^3)$, then according to \cite[Lem. 3.3, p. 264]{AMITA}, $\zeta_3$ is a norm of an element of $N-\lbrace 0\rbrace$ and $q^{*}=1$. We conclude that $\operatorname{rank}\,(Cl_{3}^{(\sigma)}(N))=3$, which is  absurd.      
   \item[2)]  If $ \exists i \in   \{1,2,3,4\}$ such that  \(q_{i} \not \equiv -1  \pmod 9\), then $D=3^{e}q_1^{f_1}q_2^{f_2}q_3^{f_3}q_4^{f_4}\equiv  \pm 3^{e}, \pm 4\times 3^{e} \ or \ \pm 7 \times 3^{e} \pmod 9 $,  since $D \equiv \pm 1 \pmod 9 $, then $e=0$, so the possible form of $D$ is $D=q_1^{f_1}q_2^{f_2}q_3^{f_3}q_4^{f_4} \equiv \pm 1 \pmod 9 $,  which is the $13^{th}$ form of $D$ in Theorem \ref{thm:Rank2}.    
   \end{itemize}
    
 %%%%%_______________
 
 \paragraph*{Species 1:} If \(D \not \equiv \pm1 \pmod  9\):

For each $i \in   \{1,2,3,4\}$, $q_{i}$ is inert in $K$ because $ q_{i} \equiv 2 \pmod 3$, and $q_{i}$ is ramified in $L$. As $D \not \equiv \pm1 \pmod  9$, $3$ is   ramified in $L$, so $\lambda$ is  ramified in $N/K$. Then $t=5$. We conclude that $\operatorname{rank}\,(Cl_{3}^{(\sigma)}(N))=3$ or $4$, which is  absurd. 

\paragraph{}
Finally, if the integer $D$ takes one of the forms given in  Theorem \ref{thm:Rank2},
 then $\operatorname{rank}\,(Cl_{3}^{(\sigma)}(N))=2$ according to the proof of Theorem \ref{thm:Ismaili1} and \ref{thm:Ismaili2}.
%------------------------------------------------------------------------

%-----------------------------------------------------------------------

\section{Conclusion}
\label{s:Conclusion}

\noindent
We have characterized all Kummer extensions $N/K$,
which possess a relative $3$-genus field $N^{\ast}$
with elementary bicyclic Galois group $\operatorname{Gal}(N^{\ast}/N)$,
by precisely thirteen forms of radicands $D$ of $N=\mathbb{Q}(\sqrt[3]{D},\zeta_3)$
with up to four prime divisors,
at most two among them congruent to $1$ modulo $3$.
However, the genus group $\operatorname{Gal}(N^{\ast}/N)$
is a common property of the whole multiplet $N_1,\ldots,N_m$
of $m$ fields sharing the same conductor $f$,
as opposed to the property of being an M0-field
\cite{AMI},
which sensibly depends on the particular shape of the radicand $D$. 
In our numerical examples, tables
\ref{tb1}, \ref{tb2}, \ref{tb3}, \ref{tb4}, and \ref{tb5},
$\operatorname{rank}\,(Cl_{3}(N))=\operatorname{rank}\,(Cl_{3}^{(\sigma)}(N))$
or $ 1+ \operatorname{rank}\,(Cl_{3}^{(\sigma)}(N))$
or $ 2+ \operatorname{rank}\,(Cl_{3}^{(\sigma)}(N))$.
Consequently, there arises the question
whether the ambiguous $3$-class group $Cl_{3}^{(\sigma)}(N)$
is a good approximation of the complete $3$-class group $Cl_{3}(N)$
or the difference
$\operatorname{rank}\,(Cl_{3}(N))-\operatorname{rank}\,(Cl_{3}^{(\sigma)}(N))$
my be arbitrarily large,
and whether the types $\alpha$ or $\gamma$ impact on
$3$-class groups $Cl_{3}(L)$ and $Cl_{3}(N)$
for being elementary or non-elementary groups.
 
%-----------------------------------------------------------------------

\section{Acknowledgements}
\label{s:Thanks}

\noindent
The fourth author acknowledges that his research was supported by the
Austrian Science Fund (FWF): projects J0497-PHY, P26008-N25,
and by the Research Executive Agency of the European Union (EUREA):
project Horizon Europe 2021--2027.
 
%-----------------------------------------------------------------------

\section{Data availability statement}
\label{s:Data}

\noindent
Data underlying our paper can be obtained from the authors upon request.
 
%-----------------------------------------------------------------------

\section{Conflict of interests statement}
\label{s:Conflict}

\noindent
The authors declare that no conflict of interest arises by publication of their paper.

%-----------------------------------------------------------------------

%-----------------------------------------------------------------------
\begin{quote}

Siham AOUISSI  \\
Ecole Normale Supérieure of Moulay Ismail University (ENS-UMI), \\
Algebraic theories and applications research team (ATA),\\
Address: ENS, B.P. 3104, Toulal, Meknès, Morocco, \\
 E-mail: aouissi.siham@gmail.com. \\

Abdelmalek AZIZI, Moulay Chrif ISMAILI \\
Department of Mathematics and Computer Sciences, \\
Mohammed first University, \\
60000 Oujda - Morocco, \\
abdelmalekazizi@yahoo.fr, mcismaili@yahoo.fr. \\

Daniel C. MAYER \\
Austrian Science Fund, \\
Naglergasse 53, 8010 Graz, Austria, \\
algebraic.number.theory@algebra.at \\
URL: http://www.algebra.at. \\

Mohamed TALBI \\
Regional Center of Professions of Education and Training, \\
60000 Oujda - Morocco, \\
ksirat1971@gmail.com.

\end{quote}
%--------------------------------------------------------------------------------

\end{document}